\documentclass[envcountsect,envcountsame]{svmult}

\usepackage[leqno]{amsmath}
\usepackage{latexsym,amsfonts,amssymb}
\usepackage[all]{xy} \SelectTips{eu}{}
\usepackage{makeidx}
\usepackage{multicol}
\usepackage{type1cm} 
\usepackage[bottom]{footmisc}
\usepackage{hyperref}

\newcommand{\rGpd}[1]{\operatorname{rel-Gpd}_R#1}
\newcommand{\GPExt}[4][GP]{\operatorname{Ext}_{\rm #1}^{#2}(#3,#4)}

\renewcommand{\ge}{\geqslant}
\renewcommand{\le}{\leqslant}

\newcommand{\setof}[2]{\{\,#1 \mid #2\,\}}
\newcommand{\ZZ}{\mathbb{Z}}
\newcommand{\NN}{\mathbb{N}}
 \newcommand{\iinZ}{{i\in\ZZ}}
 \newcommand{\f}{\varphi}
 \newcommand{\m}{\mathfrak{m}}
 \newcommand{\n}{\mathfrak{n}}
 \newcommand{\p}{\mathfrak{p}}
 \newcommand{\is}{\cong}
 \newcommand{\lra}{\longrightarrow}
 \newcommand{\xra}[2][]{\xrightarrow[#1]{\;#2\;}}
 \newcommand{\dra}[2][]{\xra{\dif[#1]{#2}}}
 \newcommand{\poly}[2][k]{#1[#2]}
 \newcommand{\pows}[2][k]{#1[\mspace{-2.3mu}[#2]\mspace{-2.3mu}]}
 \newcommand{\Rm}{(R,\m)}
 \newcommand{\Rmk}{(R,\m,k)}
 
 \newcommand{\Rhat}{\widehat{R}}
 \newcommand{\Shat}{\widehat{S}}
 \newcommand{\mapdef}[4][\rightarrow]{\nobreak{#2\colon #3 #1 #4}}
 \newcommand{\dmapdef}[4][\lra]{\nobreak{#2\colon #3\:#1\:#4}}
 \renewcommand{\Im}[1]{\nobreak{\operatorname{Im}#1}}
 \newcommand{\Ker}[1]{\nobreak{\operatorname{Ker}#1}}
 \newcommand{\Coker}[1]{\nobreak{\operatorname{Coker}#1}}
 \newcommand{\dif}[2][]{{\partial}^{#2}_{#1}}
 \renewcommand{\H}[2][]{\operatorname{H}_{#1}(#2)}
 \newcommand{\HH}[2][]{\operatorname{H}^{#1}(#2)}
 \newcommand{\bas}[3][R]{\mu_{#1}^{#2}(#3)}
 \newcommand{\SpecR}{\operatorname{Spec}R}
 \newcommand{\dptR}{\operatorname{depth}R}
 \newcommand{\dimR}{\operatorname{dim}R}
 \newcommand{\edim}[1]{\operatorname{edim}#1}
 \newcommand{\Ann}[2][R]{\operatorname{Ann}_{#1}#2}
 \newcommand{\Supp}[2][R]{\operatorname{Supp}_{#1}#2}
 \newcommand{\wdt}[2][R]{\operatorname{width}_{#1}#2}
 \newcommand{\dpt}[2][R]{\operatorname{depth}_{#1}#2}
 \newcommand{\grd}[2][R]{\operatorname{grade}_{#1}#2}
 \newcommand{\fd}[2][R]{\operatorname{fd}_{#1}#2}
 \newcommand{\id}[2][R]{\operatorname{id}_{#1}#2}
 \newcommand{\pd}[2][R]{\operatorname{pd}_{#1}#2}
 \newcommand{\Gdim}[2][R]{\operatorname{G-dim}_{#1}#2}
 \newcommand{\Gfd}[2][R]{\operatorname{Gfd}_{#1}#2}
 \newcommand{\Gid}[2][R]{\operatorname{Gid}_{#1}#2}
 \newcommand{\Gpd}[2][R]{\operatorname{Gpd}_{#1}#2}
 \newcommand{\Hom}[3][R]{\operatorname{Hom}_{#1}(#2,#3)}
 \newcommand{\RHom}[3][R]{\operatorname{\mathbf{R}Hom}_{#1}(#2,#3)}
 \newcommand{\Ext}[4][R]{\operatorname{Ext}_{#1}^{#2}(#3,#4)}

 \newcommand{\tp}[3][R]{\nobreak{#2\otimes_{#1}#3}}
 \newcommand{\tpp}[3][R]{(\tp[#1]{#2}{#3})}
 
 \newcommand{\Ltp}[3][R]{\nobreak{#2\otimes_{#1}^{\mathbf{L}}#3}}

 \newcommand{\Tor}[4][R]{\operatorname{Tor}^{#1}_{#2}(#3,#4)}
 \newcommand{\Catfont}[1]{\mathcal{#1}}
 \newcommand{\Cat}[2]{{\Catfont{#2}}(#1)}
 \newcommand{\Du}[1][R]{\Cat{#1}{D}}
 \newcommand{\R}[1][R]{\Cat{#1}{R}}
 \newcommand{\A}[1][R]{\Cat{#1}{A}}
 \newcommand{\AD}[1][D]{\Catfont{A}_{#1}(R)}
 \newcommand{\BD}[1][D]{\Catfont{B}_{#1}(R)}
 \newcommand{\Au}[1][R]{\Cat{#1}{\widehat{A}}}
 \newcommand{\B}[1][R]{\Cat{#1}{B}}
 \newcommand{\Bu}[1][R]{\Cat{#1}{\widehat{B}}}

\newcommand{\thmref}[2][Theorem~]{#1\ref{thm:#2}}
\newcommand{\corref}[2][Corollary~]{#1\ref{cor:#2}}
\newcommand{\prpref}[2][Proposition~]{#1\ref{prp:#2}}
\newcommand{\lemref}[2][Lemma~]{#1\ref{lem:#2}}

\newcommand{\prbref}[2][Problem~]{#1\ref{prb:#2}}
\newcommand{\qstref}[2][Question~]{#1\ref{qst:#2}}
\newcommand{\dfnref}[2][Definition~]{#1\ref{dfn:#2}}
\newcommand{\exaref}[2][Example~]{#1\ref{exa:#2}}
\newcommand{\rmkref}[2][Remark~]{#1\ref{rmk:#2}}

\newcommand{\secref}[2][Section~]{#1\ref{sec:#2}}
\newlength{\mylabelsep}
\newcommand{\cH}{\Catfont{H}}
\newcommand{\cx}[1]{\boldsymbol{#1}}

\spnewtheorem{construction}[theorem]{Construction}{\it}{\rm}
\spnewtheorem{prb}[theorem]{Problem}{\it}{\rm}
\spnewtheorem*{resa}{Theorem~A}{\it}{\rm}
\spnewtheorem*{resb}{Theorem~B}{\it}{\rm}
\spnewtheorem*{resc}{Theorem~C}{\it}{\rm}
\spnewtheorem*{resd}{Theorem~D}{\it}{\rm}
\spnewtheorem*{rese}{Theorem~E}{\it}{\rm}
\spnewtheorem*{resf}{Theorem~F}{\it}{\rm}
\spnewtheorem*{resg}{Theorem~G}{\it}{\rm}
\spnewtheorem*{resh}{Theorem~H}{\it}{\rm}
\spnewtheorem*{metatheorem}{Meta Question}{\it}{\rm}
\spnewtheorem{dfn}[theorem]{Definition}{\bf}{\rm}

\makeindex             

\begin{document}
\smartqed \title*{Beyond Totally Reflexive Modules and~Back}

\subtitle{A Survey on Gorenstein Dimensions}

\titlerunning{Beyond Totally Reflexive Modules and Back}

\toctitle{Beyond Totally Reflexive Modules and Back}

\author{Lars Winther Christensen, Hans-Bj\o rn Foxby, and Henrik Holm}

\authorrunning{L.\,W.~Christensen\and H.-B.~Foxby\and H.~Holm}

\institute{L.\,W.\ Christensen \at Department of Mathematics and
  Statistics, Texas Tech University, Mail Stop 1042, Lubbock,
  TX~79409, U.S.A. %
  \email{lars.w.christensen@ttu.edu}\and %
  H.-B.~Foxby \at Department of Mathematical Sciences, University of
  Copenhagen, Universitets\-parken~5, DK-2100 K\o benhavn \O, Denmark
  \email{foxby@math.ku.dk}\and %
  H.\ Holm \at Department of Basic Sciences and Environment,
  University of Copenhagen, Thorvaldsensvej 40, DK-1871 Frederiksberg
  C, Denmark \email{hholm@life.ku.dk} }

%
%

\maketitle

\vspace*{-3.2\baselineskip}

\abstract{Starting from the notion of totally reflexive modules, we
  survey the theory of Gorenstein homological dimensions for modules
  over commutative rings. The account includes the theory's
  connections with relative homological algebra and with studies of
  local ring homomorphisms. It ends close to the starting point: with
  a characterization of Gorenstein rings in terms of total acyclicity
  of complexes.  \keywords{Auslander categories, Bass class, Chouinard
    formula, Frobenius endomorphism, G-dimension, Gorenstein
    dimension, Gorenstein flat cover, Gorenstein injective
    preenvelope, Gorenstein projective precover, Gorenstein ring,
    quasi-Cohen--Macaulay homomorphism, quasi-Gorenstein homomorphism,
    totally acyclic complex, totally
    reflexive module.\\[\baselineskip]
    \textbf{Mathematics Subject Classification (2000):}\newline
    13--02, 13B10, 13D05, 13H10, 18G25.\\[.2\baselineskip]} }

\setcounter{minitocdepth}{2} \dominitoc


\section*{Introduction}

An important motivation for the study of homological dimensions dates
back to 1956, when Auslander and Buchsbaum~\cite{MAsDBc56} and
Serre~\cite{JPS56} proved:

  \begin{resa}
    Let $R$ be a commutative Noetherian local ring with residue
    field~$k$. Then the following conditions are equivalent.
    \begin{description}[$(iii)$]\setlength{\labelsep}{\mylabelsep}
    \item[{\hfill $(i)$}] $R$ is regular. %
      \index{regular local ring}
    \item[\hfill $(ii)$] $k$ has finite projective dimension.
    \item[\hfill $(iii)$] Every $R$-module has finite projective
      dimension.
    \end{description}
  \end{resa}

\noindent
This result opened to the solution of two long-standing conjectures of
Krull.  Moreover, it introduced the theme that finiteness of a
homological dimension for all modules characterizes rings with special
properties.  Later work has shown that modules of finite projective
dimension over a general ring share many properties with modules over
a regular ring. This is an incitement to study homological dimensions
of individual modules.

In line with these ideas, Auslander and Bridger \cite{MAsMBr69}
introduced in 1969~the G-dimension. It is a homological dimension for
finitely generated modules over a Noetherian ring, and it gives a
characterization of Gorenstein local rings (\thmref{Gor-gdim}), which
is similar to Theorem~A. Indeed, $R$ is Gorenstein if $k$ has finite
G-dimension, and only if every finitely generated $R$-module has
finite G-dimension.

In the early 1990s, the G-dimension was extended beyond the realm of
finitely generated modules over a Noetherian ring. This was done by
Enochs and Jenda who introduced the notion of Gorenstein projective
modules~\cite{EEnOJn95b}.  With the Gorenstein projective dimension at
hand, a perfect parallel to Theorem~A becomes available
(\thmref{Gor-Gpd}). Subsequent work has shown that modules of finite
Gorenstein projective dimension over a general ring share many
properties with modules over a Gorenstein ring.

\subsubsection*{Classical Homological Algebra as Precedent}

The notions of injective dimension and flat dimension for modules also
have Gorenstein counterparts. It was Enochs and Jenda who introduced
Gorenstein injective modules~\cite{EEnOJn95b} and, in collaboration
with Torrecillas, Gorenstein flat modules~\cite{EJT-93}.  The study of
Gorenstein dimensions is often called \emph{Gorenstein homological
  algebra}; it has taken directions from the following:

\begin{metatheorem}
  Given a result in classical homological algebra, does it have a
  counterpart in Gorenstein homological algebra?
\end{metatheorem}

\noindent
To make this concrete, we review some classical results on homological
dimensions and point to their Gorenstein counterparts within the main
text.  In the balance of this introduction, $R$ is assumed to be a
commutative Noetherian local ring with maximal ideal $\m$ and residue
field $k=R/\m$.

\paragraph{\sc Depth and Finitely Generated Modules}

The projective dimension of a finitely generated $R$-module is closely
related to its depth. This is captured by the
\emph{Auslander--Buchsbaum Formula}~\cite{MAsDBc57}: %
\index{Auslander--Buchsbaum Formula}

\begin{resb}
  For every finitely generated $R$-module $M$ of finite projective
  dimension there is an equality $\pd{M} = \dptR - \dpt{M}.$
\end{resb}

\noindent
The Gorenstein counterpart (\thmref{ABr}) actually strengthens the
classical result; this is a recurring theme in Gorenstein homological
algebra.

The injective dimension of a non-zero finitely generated $R$-module is
either infinite or it takes a fixed value:

\begin{resc}
  For every non-zero finitely generated $R$-module $M$ of finite
  injective dimension there is an equality $\id{M} = \dptR.$ %
  \index{Bass Formula}
\end{resc}

\noindent
This result of Bass~\cite{HBs63} has its Gorenstein counterpart in
\thmref{Bass_formula}.

\paragraph{\sc Characterizations of Cohen--Macaulay
  Rings}

Existence of special modules of finite homological dimension
characterizes Cohen--Macaulay rings. The equivalence of $(i)$ and
$(iii)$ in the theorem below is still referred to as \emph{Bass'
  conjecture,} %
\index{Bass' conjecture} even though it was settled more than 20 years
ago. Indeed, Peskine and Szpiro proved in \cite{CPsLSz73} that it
follows from the New Intersection Theorem, which they proved
\emph{ibid.}\ for equicharacteristic rings. In 1987
Roberts~\cite{PRb87} settled the New Intersection Theorem completely.
\begin{resd}
  The following conditions on $R$ are equivalent.
  \begin{description}[$(iii)$]\setlength{\labelsep}{\mylabelsep}
  \item[\hfill $(i)$] $R$ is Cohen--Macaulay. %
    \index{Cohen--Macaulay local ring}
  \item[\hfill $(ii)$] There is a non-zero $R$-module of finite length
    and finite projective~dim.
  \item[\hfill $(iii)$] There is a non-zero finitely generated
    $R$-module of finite injective dim.
  \end{description}
\end{resd}

\noindent
A Gorenstein counterpart to this characterization is yet to be
established; see \qstref[Questions~]{gdim-cm} and~\qstref[]{gid-cm}.

Gorenstein rings are also characterized by existence of special
modules of finite homological dimension. The equivalence of $(i)$ and
$(ii)$ below is due to Peskine and Szpiro~\cite{CPsLSz73}. The
equivalence of $(i)$ and $(iii)$ was conjectured by
Vasconcelos~\cite{dtmc} and proved by Foxby~\cite{HBF77b}. The
Gorenstein counterparts are given in
\thmref[Theorems~]{localGorenstein} and \thmref[]{Gor-Gfd}; see also
\qstref{GHBF}.

\begin{rese}
  The following conditions on $R$ are equivalent.
  \begin{description}[$(iii)$]\setlength{\labelsep}{\mylabelsep}
  \item[\hfill $(i)$] $R$ is Gorenstein. %
    \index{Gorenstein ring!local}
  \item[\hfill $(ii)$] There is a non-zero cyclic $R$-module of finite
    injective dimension.
  \item[\hfill $(iii)$] There is a non-zero finitely generated
    $R$-module of finite projective dimension and finite injective
    dimension.
  \end{description}
\end{rese}

\paragraph{\sc Local Ring Homomorphisms}

The Frobenius endomorphism detects regularity of a local ring of
positive prime characteristic. The next theorem collects results of
Avramov, Iyengar, and Miller~\cite{AIM-06}, Kunz~\cite{EKn69}, and
Rodicio~\cite{AGR88}. The counterparts in Gorenstein homological
algebra to these results are given in \thmref[Theorems~]{FrobGfd} and
\thmref[]{FrobGid}.

\begin{resf}
  Let $R$ be of positive prime characteristic, and let $\phi$ denote
  its Frobenius endomorphism. Then the following conditions are
  equivalent. %
  \index{Frobenius endomorphism}
  \begin{description}[$(iii)$]\setlength{\labelsep}{\mylabelsep}
  \item[\hfill $(i)$] $R$ is regular. %
    \index{regular local ring}
  \item[\hfill $(ii)$] $R$ has finite flat dimension as an $R$-module
    via $\phi^n$ for some $n \ge 1$.
  \item[\hfill $(iii)$] $R$ is flat as an $R$-module via $\phi^n$ for
    every $n \ge 1$.
  \item[\hfill $(iv)$] $R$ has finite injective dimension as an
    $R$-module via $\phi^n$ for some $n \ge 1$.
  \item[\hfill $(v)$] $R$ has injective dimension equal to $\dimR$ as
    an $R$-module via $\phi^n$ for every~\mbox{$n \ge 1$}.
  \end{description}
\end{resf}

Let $(S,\n)$ be yet a commutative Noetherian local ring. A ring
homomorphism $\mapdef{\f}{R}{S}$ is called \emph{local} %
\index{local homomorphism} if there is an inclusion \mbox{$\f(\m)
  \subseteq \n$}. A classical chapter of local algebra, initiated by
Grothendieck, studies transfer of ring theoretic properties along such
homomorphisms. If $\f$ is flat, then it is called
\emph{Cohen--Macaulay} or \emph{Gorenstein} if its closed fiber $S/\m
S$ is, respectively, a Cohen--Macaulay ring or a Gorenstein
ring. These definitions have been extended to homomorphisms of finite
flat dimension. The theorem below collects results of Avramov and
Foxby from \cite{LLAHBF92} and \cite{LLAHBF98}; the Gorenstein
counterparts are given in \thmref[Theorems~]{GorAscDescFtGFD} and
\thmref[]{CMAscDsc}.

\begin{resg}
  Let $\mapdef{\f}{R}{S}$ be a local homomorphism and assume that $S$
  has finite flat dimension as an $R$-module via $\f$. Then the
  following hold:
  \begin{description}[$(a)$]\setlength{\labelsep}{\mylabelsep}
  \item[\hfill \rm (a)] $S$ is Cohen--Macaulay if and only if $R$ and
    $\f$ are Cohen--Macaulay. %
    \index{local homomorphism!Cohen--Macaulay}
  \item[\hfill \rm (b)] $S$ is Gorenstein if and only if $R$ and $\f$
    are Gorenstein. %
    \index{local homomorphism!Gorenstein}
  \end{description}
\end{resg}

\paragraph{\sc Vanishing of Cohomology}

The projective dimension of a module $M$ is at most $n$ if and only if
the absolute cohomology functor $\Ext[]{n+1}{M}{-}$
vanishes. Similarly (\thmref{GPExt}), $M$ has Gorenstein projective
dimension at most $n$ if and only if the relative cohomology functor
$\Ext[\rm GP]{n+1}{M}{-}$ vanishes. Unfortunately, the similarity
between the two situations does not run too deep. We give a couple of
examples:

The absolute Ext is balanced, that is, it can be computed from a
projective resolution of $M$ or from an injective resolution of the
second argument. In general, however, the only known way to compute
the relative Ext is from a (so-called) proper Gorenstein projective
resolution of $M$.

Secondly, if $M$ is finitely generated, then the absolute Ext commutes
with localization, but the relative Ext is not known to do so, unless
$M$ has finite Gorenstein projective dimension.

Such considerations motivate the search for an alternative
characterization of modules of finite Gorenstein projective dimension,
and this has been a driving force behind much research on Gorenstein
dimensions within the past 15 years. What follows is a brief review.

\paragraph{\sc Equivalence of Module Categories}

For a finitely generated $R$-module, Foxby \cite{HBF94} gave a
``resolution-free'' criterion for finiteness of the Gorenstein
projective dimension; that is, one that does not involve construction
of a Gorenstein projective resolution. This result from 1994 is
\thmref{Rmod}. In 1996, Enochs, Jenda, and Xu~\cite{EJX-96b} extended
Foxby's criterion to non-finitely generated $R$-modules, provided that
$R$ is Cohen--Macaulay with a dualizing module $D$. Their work is
related to a 1972 generalization by Foxby~\cite{HBF72} of a theorem of
Sharp~\cite{RYS72}. Foxby's version~reads:

\begin{resh}
  Let $R$ be Cohen--Macaulay with a dualizing module $D$.  Then the
  horizontal arrows below are equivalences of categories of
  $R$-modules. \vspace*{-.8\baselineskip}
  \begin{petit}
    $$\xymatrix@C=8em{ \AD \ar@<0.75ex>[r]^{\tp{D}{-}} 
      & \BD \ar@<.75ex>[l]^{\Hom D-}\\
      \{ A \mid \pd{A} \text{ is finite} \}
      \ar@<0.75ex>[r]^{\tp{D}{-}} \ar@{^(->}[u] & \{ B \mid \id{B}
      \text{ is finite} \} \ar@<.75ex>[l]^{\Hom D-}\ar@{^(->}[u] }$$
  \end{petit}
  \index{Auslander class}%
  \index{Bass class}%
  \index{dualizing complex}%
  \index{Cohen--Macaulay local ring}
\end{resh}
\noindent
Here $\AD$ is the \emph{Auslander class} (\dfnref{A}) with respect to
$D$ and $\BD$ is the \emph{Bass class} (\dfnref{B}). What Enochs,
Jenda, and Xu prove in \cite{EJX-96b} is that the $R$-modules of
finite Gorenstein projective dimension are exactly those in the $\AD$,
and the modules in $\BD$ are exactly those of finite Gorenstein
injective dimension. Thus, the upper level equivalence in Theorem~H is
the Gorenstein counterpart of the lower level equivalence.

A commutative Noetherian ring has a dualizing complex $\cx{D}$ if and
only if it is a homomorphic image of a Gorenstein ring of finite Krull
dimension; see Kawasaki~\cite{TKw02}. For such rings, a result similar
to Theorem~H was proved by Avramov and Foxby~\cite{LLAHBF97} in 1997.
An interpretation in terms of Gorenstein dimensions
(\thmref[Theorems~]{Amod}~and~\thmref[]{Bmod}) of the objects in
$\AD[\cx{D}]$ and $\BD[\cx{D}]$ was established by Christensen,
Frankild, and Holm~\cite{CFH-06} in 2006. Testimony to the utility of
these results is the frequent occurrence---e.g.\ in
\thmref[Theorems~]{GIcolimit}, \thmref[]{GFprod}, \thmref[]{noname},
\thmref[]{HHlPJr}, \thmref[]{FrobGid}, \thmref[]{gdim-Gfd},
\thmref[]{GorAD}, and \thmref[]{qGorAD}---of the assumption that the
ground ring is a homomorphic image of a Gorenstein ring of finite
Krull~dimension. Recall that every complete local ring satisfies this
assumption.

Recent results, \thmref[Theorems~]{GPcompletion} and
\thmref[]{GFcompletion}, by Esmkhani and Tousi~\cite{MAEMTs07} and
\thmref{GidCoBase} by Christensen and Sather-Wagstaff~\cite{LWCSSWc}
combine with \thmref[Theorems~]{Amod}~and~\thmref[]{Bmod} to provide
resolution-free criteria for finiteness of Gorenstein dimensions over
general local rings; see \rmkref[Remarks~]{resfree1} and
\rmkref[]{Gid}.

\subsubsection*{Scope and Organization}

A survey of this modest length is a portrait painted with broad pen
strokes. Inevitably, many details are omitted, and some generality has
been traded in for simplicity. We have chosen to focus on modules over
commutative, and often Noetherian, rings. Much of Gorenstein
homological algebra, though, works flawlessly over non-commutative
rings, and there are statements in this survey about Noetherian rings
that remain valid for coherent rings. Furthermore, most statements
about modules remain valid for complexes of modules. The reader will
have to consult the references to qualify these claims.

In most sections, the opening paragraph introduces the main references
on the topic. We strive to cite the strongest results available and,
outside of this introduction, we do not attempt to trace the history
of individual results.  In notes, placed at the end of sections, we
give pointers to the literature on directions of research---often new
ones---that are not included in the survey. Even within the scope of
this paper, there are open ends, and more than a dozen questions and
problems are found throughout the text.

From this point on, $R$ denotes a commutative ring. Any extra
assumptions on $R$ are explicitly stated. We say that $R$ is
\emph{local} %
\index{local ring} if it is Noetherian and has a unique maximal
ideal. We use the shorthand $\Rmk$ for a local ring $R$ with maximal
ideal $\m$ and residue field $k=R/\m$.

\subsubsection*{Acknowledgments}

It is a pleasure to thank Nanqing Ding, his colleagues, and the
students who attended the 2008 Summer School at Nanjing
University. The choice of material and the organization of this survey
is strongly influenced by the summer school participants' reactions to
lectures by L.W.C. We also thank Edgar Enochs, Srikanth Iyengar, and
Peter J{\o}rgensen for helpful discussions on their past and current
work related to Gorenstein homological algebra. In the preparation of
this final version, we have been helped and challenged by a series of
detailed and thoughtful comments from our colleagues. We thank Sean
Sather-Wagstaff and everyone else who gave their comments on earlier
versions of the paper: Driss Bennis, Jesse Burke, Nanqing Ding, Najib
Mahdou, Ryo Takahashi, and Oana Veliche. We also acknowledge the
referee's careful reading of the manuscript and many helpful remarks.

This paper was started and finalized during visits by L.W.C.\ to
Copenhagen; the hospitality and support of the University of
Copenhagen is gratefully acknowledged.


\section{G-dimension of Finitely Generated Modules}
\label{sec:gdim}

The topic of this section is Auslander and Bridger's notion of
\mbox{G-dimension} for finitely generated modules over a Noetherian
ring.  The notes \cite{MAs67} from a seminar by Auslander outline the
theory of G-dimension over commutative Noetherian rings. In
\cite{MAsMBr69} Auslander and Bridger treat the G-dimension within a
more abstract framework. Later expositions are given by
Christensen~\cite{lnm} and by Ma\c{s}ek~\cite{VMs00}.

A \emph{complex} %
\index{complex} $\cx{M}$ of modules is (in homological notation) an
infinite sequence of homomorphisms of $R$-modules
$$\cx{M}= \cdots \dra[i+1]{} M_{i} \dra[i]{} M_{i-1} \dra[i-1]{} \cdots$$
such that $\dif[i]{}\dif[i+1]{}=0$ for every $i$ in $\ZZ$. The $i\,$th
\emph{homology module} of $\cx{M}$ is $\H[i]{\cx{M}} =
\Ker{\dif[i]{}}/\Im{\dif[i+1]{}}$. We call $\cx{M}$ \emph{acyclic} %
\index{acyclic complex} if $\H[i]{\cx{M}} =0$ for all $i\in\ZZ$.

\begin{lemma}
  \label{lem:taL}
  Let $\cx{L}$ be an acyclic complex of finitely generated projective
  $R$-modules. The following conditions on $\cx{L}$ are equivalent.
  \begin{description}[$(iii)$]\setlength{\labelsep}{\mylabelsep}
  \item[\hfill $(i)$] The complex $\Hom{\cx{L}}{R}$ is acyclic.
  \item[\hfill $(ii)$] The complex $\Hom{\cx{L}}{F}$ is acyclic for
    every flat $R$-module $F$.
  \item[\hfill $(iii)$] The complex $\tp{E}{\cx{L}}$ is acyclic for
    every injective $R$-module $E$.
  \end{description}
\end{lemma}

\begin{petit}
  \begin{proof}
    The Lemma is proved in \cite{lnm}, but here is a cleaner argument:
    Let $F$ be a flat module and $E$ be an injective module. As
    $\cx{L}$ consists of finitely generated projective modules, there
    is an isomorphism of complexes
    $$\Hom{\Hom{\cx{L}}{F}}{E} \is \tp{\Hom{F}{E}}{\cx{L}}\,.$$
    It follows from this isomorphism, applied to $F=R$, that $(i)$
    implies $(iii)$. Applied to a faithfully injective module $E$, it
    shows that $(iii)$ implies $(ii)$, as $\Hom{F}{E}$ is an injective
    module. It is evident that $(ii)$ implies $(i)$.\qed
  \end{proof}
\end{petit}

The following nomenclature is due to Avramov and Martsinkovsky
\cite{LLAAMr02}; \lemref{totref} clarifies the rationale behind it.

\begin{dfn}
  \label{dfn:taL}
  A complex $\cx{L}$ that satisfies the conditions in \lemref{taL} is
  called \emph{totally acyclic}. %
  \index{totally acyclic complex} An $R$-module $M$ is called
  \emph{totally reflexive} %
  \index{totally reflexive!module} if there exists a totally acyclic
  complex $\cx{L}$ such that $M$ is isomorphic to $\Coker{(L_1 \to
    L_0)}$.
\end{dfn}

\noindent
Note that every finitely generated projective module $L$ is totally
reflexive; indeed, the complex \mbox{$0 \to L \xra{=} L \to 0$}, with
$L$ in homological degrees $0$ and $-1$, is totally acyclic.

\begin{example}
  \label{exa:uv}
  If there exist elements $x$ and $y$ in $R$ such that
  \mbox{$\Ann{(x)}=(y)$} and \mbox{$\Ann{(y)}=(x)$}, then the complex
  \begin{displaymath}
    \cdots \xra{x} R \xra{y} R  \xra{x} R \xra{y} \cdots
  \end{displaymath}
  is totally acyclic. Thus, $(x)$ and $(y)$ are totally reflexive
  $R$-modules. For instance, if $X$ and $Y$ are non-zero non-units in
  an integral domain $D$, then their residue classes $x$ and $y$ in
  $R=D/(XY)$ generate totally reflexive $R$-modules.
\end{example}

\noindent
An elementary construction of rings of this kind---\exaref{taL}
below---shows that non-projective totally reflexive modules may exist
over a variety of rings; see also \prbref{LLAAMr}.

\begin{example}
  \label{exa:taL}
  Let $S$ be a commutative ring, and let $m>1$ be an integer. Set
  $R=\poly[S]{X}/(X^m)$, and denote by $x$ the residue class of $X$ in
  $R$. Then for every integer $n$ between $1$ and $m-1$, the module
  $(x^n)$ is totally reflexive.
\end{example}

From \lemref{taL} it is straightforward to deduce:

\begin{proposition}
  \label{prp:1}
  Let $S$ be an $R$-algebra of finite flat dimension. For every
  totally reflexive $R$-module $G$, the module $\tp{S}{G}$ is totally
  reflexive over $S$.
\end{proposition}

\noindent
\prpref{1} applies to \mbox{$S=R/(\boldsymbol{x})$}, where
$\boldsymbol{x}$ is an $R$-regular sequence. If $(R,\m)$ is local,
then it also applies to the $\m$-adic completion~\mbox{$S=\Rhat$.}

\subsubsection*{Noetherian Rings}

Recall that a finitely generated $R$-module $M$ is called
\emph{reflexive} %
\index{reflexive module} if the canonical map from $M$ to
$\Hom{\Hom{M}{R}}{R}$ is an isomorphism. The following
characterization of totally reflexive modules goes back to
\cite[4.11]{MAsMBr69}.

\begin{lemma}
  \label{lem:totref}
  Let $R$ be Noetherian. A finitely generated $R$-module $G$ is
  totally reflexive if and only if it is reflexive and for every $i
  \ge 1$ one has
  $$\Ext{i}{G}{R} = 0 = \Ext{i}{\Hom{G}{R}}{R}.$$
\end{lemma}

\begin{dfn}
  An (augmented) \emph{G-resolution} of a finitely generated module
  $M$ is an exact sequence $\cdots \to G_i \to G_{i-1} \to \cdots \to
  G_0 \to M \to 0,$ where each module $G_i$ is totally reflexive.
\end{dfn}

\noindent
Note that if $R$ is Noetherian, then every finitely generated
$R$-module has a G-resolution, indeed it has a resolution by finitely
generated free modules.

\begin{dfn}
  \label{dfn:G-dim}
  Let $R$ be Noetherian. For a finitely generated $R$-module \mbox{$M
    \ne 0$} the \emph{G-dimension}, %
  \index{G-dimension} denoted by $\Gdim{M}$, is the least integer $n
  \ge 0$ such that there exists a G-resolution of $M$ with $G_i=0$ for
  all $i > n$. If no such $n$ exists, then $\Gdim{M}$ is infinite. By
  convention, set \mbox{$\Gdim{0} = -\infty$}.
\end{dfn}

\noindent
The `G' in the definition above is short for Gorenstein.

In \cite[ch.~3]{MAsMBr69} one finds the next theorem and its
corollary; see also \cite[1.2.7]{lnm}.

\pagebreak[4]

\begin{theorem}
  \label{thm:gdim}
  Let $R$ be Noetherian and $M$ be a finitely generated $R$-module of
  finite G-dimension. For every $n\ge 0$ the next conditions are
  equivalent.
  \begin{description}[$(iii)$]\setlength{\labelsep}{\mylabelsep}
  \item[\hfill $(i)$] $\Gdim{M} \le n$.
  \item[\hfill $(ii)$] $\Ext{i}{M}{R}=0$ for all $i > n$.
  \item[\hfill $(iii)$] $\Ext{i}{M}{N}=0$ for all $i > n$ and all
    $R$-modules $N$ with $\fd{N}$ finite.
  \item[\hfill $(iv)$] In every augmented G-resolution $$\cdots \to
    G_i \to G_{i-1} \to \cdots \to G_0 \to M \to 0$$ the module
    $\Coker{(G_{n+1} \to G_n)}$ is totally reflexive.
  \end{description}
\end{theorem}

\begin{corollary}
  \label{cor:gdim}
  Let $R$ be Noetherian. For every finitely generated $R$-module $M$
  of finite G-dimension there is an equality
  $$\Gdim{M} = \sup\setof{\iinZ}{\Ext{i}{M}{R} \ne 0}.$$
\end{corollary}

\begin{remark}
  Examples due to Jorgensen and \c{S}ega~\cite{DAJLMS06} show that in
  \corref{gdim} one cannot avoid the \emph{a priori} condition that
  $\Gdim{M}$ is finite.
\end{remark}

\begin{remark}
  For a module $M$ as in \corref{gdim}, the small finitistic
  projective dimension of $R$ is an upper bound for $\Gdim{M}$;
  cf.~Christensen and Iyengar~\cite[3.1(a)]{LWCSIn07}.
\end{remark}

A standard argument, see \cite[3.16]{MAsMBr69} or
\cite[3.4]{LLAAMr02}, yields:

\begin{proposition}
  \label{prp:3mg}
  Let $R$ be Noetherian. If any two of the modules in an exact
  sequence \mbox{$0 \to M' \to M \to M'' \to 0$} of finitely generated
  $R$-modules have finite G-dimension, then so has the third.
\end{proposition}

The following quantitative comparison establishes the G-dimension~as a
refinement of the projective dimension for finitely generated
modules. It is easily deduced from \corref{gdim};
see~\cite[1.2.10]{lnm}.

\begin{proposition}
  Let $R$ be Noetherian. For every finitely generated $R$-mod\-ule $M$
  one has $\Gdim{M} \le \pd{M}$, and equality holds if\, $\pd{M}$
  is~finite.
\end{proposition}

By \cite[4.15]{MAsMBr69} the G-dimension of a module can be measured
locally:

\begin{proposition}
  \label{prp:gdimp}
  Let $R$ be Noetherian. For every finitely generated $R$-mod\-ule $M$
  there is an equality $\Gdim{M} =
  \sup\setof{\Gdim[R_\p]{M_\p}}{\p\in\SpecR}$.
\end{proposition}

\noindent For the projective dimension even more is known: Bass and
Murthy~\cite[4.5]{HBsMPM67} prove that if a finitely generated module
over a Noetherian ring has finite projective dimension locally, then
it has finite projective dimension globally---even if the ring has
infinite Krull dimension. A Gorenstein counterpart has recently been
established by Avramov, Iyengar, and Lipman~\cite[6.3.4]{AIL-}.

\begin{theorem}
  \label{thm:ail}
  Let $R$ be Noetherian and let $M$ be a finitely generated
  $R$-module. If $\Gdim[R_\m]{M_\m}$ is finite for every maximal ideal
  $\m$ in $R$, then $\Gdim{M}$ is finite.
\end{theorem}

Recall that a local ring is called \emph{Gorenstein}%
\index{Gorenstein ring!local}%
\index{Gorenstein ring} if it has finite self-injective dimension.  A
Noetherian ring is Gorenstein if its localization at each prime ideal
is a Gorenstein local ring, that is, $\id[R_\p]{R_\p}$ is finite for
every prime ideal $\p$ in $R$. Consequently, the self-injective
dimension of a Gorenstein ring equals its Krull dimension; that is
$\id{R} = \dimR$. The next result follows from
\mbox{\cite[4.20]{MAsMBr69}} in combination with \prpref{gdimp}.

\begin{theorem}
  \label{thm:GorG}
  Let $R$ be Noetherian and $n \ge 0$ be an integer. Then $R$ is
  Gorenstein with $\dimR \le n$ if and only if one has $\Gdim{M} \le
  n$ for every finitely generated $R$-module $M$.
\end{theorem}

\noindent
A corollary to \thmref{ail} was established by Goto~\cite{SGt82}
already in 1982; it asserts that also Gorenstein rings of infinite
Krull dimension are characterized by finiteness of G-dimension.

\begin{theorem}
  \label{thm:Goto}
  Let $R$ be Noetherian. Then $R$ is Gorenstein if and only if every
  finitely generated $R$-module has finite G-dimension.
\end{theorem}

Recall that the \emph{grade} %
\index{grade} of a finitely generated module $M$ over a Noetherian
ring $R$ can be defined as follows: $$\grd{M} =
\inf\setof{\iinZ}{\Ext{i}{M}{R} \ne 0}.$$ Foxby~\cite{HBF75b} makes
the following:

\begin{dfn}
  \label{dfn:qperf}
  Let $R$ be Noetherian. A finitely generated $R$-module $M$ is called
  \emph{quasi-perfect} %
  \index{quasi-perfect module} if it has finite G-dimension equal to
  $\grd{M}$.
\end{dfn}

The next theorem applies to \mbox{$S=R/(\boldsymbol{x})$}, where
$\boldsymbol{x}$ is an $R$-regular sequence.  Special (local) cases of
the theorem are due to Golod \cite{ESG84} and to Avramov and Foxby
\cite[7.11]{LLAHBF97}. Christensen's proof \cite[6.5]{LWC01a}
establishes the general case.

\begin{theorem}
  \label{thm:trans}
  Let $R$ be Noetherian and $S$ be a commutative Noetherian
  module-finite $R$-algebra. If $S$ is a quasi-perfect $R$-module of
  grade $g$ such that $\Ext{g}{S}{R} \is S$, then the next equality
  holds for every finitely generated $S$-module $N$, $$\Gdim{N} =
  \Gdim[S]{N} + \Gdim{S}.$$
\end{theorem}

\noindent
Note that an $S$-module has finite G-dimension over $R$ if and only if
it has finite G-dimension over $S$; see also \thmref{qGorAD}. The next
question is raised~in~\cite{LLAHBF97}; it asks if the assumption of
quasi-perfectness in \thmref{trans} is necessary.

\begin{question}
  \label{qst:trans}
  Let $R$ be Noetherian, let $S$ be a commutative Noetherian
  module-finite $R$-algebra, and let $N$ be a finitely generated
  $S$-module. If $\Gdim[S]{N}$ and $\Gdim{S}$ are finite, is then
  $\Gdim{N}$ finite?
\end{question}

\noindent
This is known as the \emph{Transitivity Question.} %
\index{G-dimension!Transitivity Question} By \cite[4.7]{LLAHBF97} and
\cite[3.15 and 6.5]{LWC01a} it has an affirmative answer if
$\pd[S]{N}$ is finite; see also \thmref{compos}.

\subsubsection*{Local Rings}
\label{sec:gdiml}

Before we proceed with results on G-dimension of modules over local
rings, we make a qualitative comparison to the projective dimension.
\thmref{trans} reveals a remarkable property of the G-dimension, one
that has almost no counterpart for the projective dimension. Here is
an example:

\begin{example}
  Let $\Rmk$ be local of positive depth. Pick a regular element $x$ in
  $\m$ and set $S = R/(x)$. Then one has $\grd{S} = 1 = \pd{S}$ and
  $\Ext{1}{S}{R} \is S$, but $\pd[S]{N}$ is infinite for every
  $S$-module $N$ such that $x$ is in $\m\Ann{N}$; see Shamash \cite[\S
  3]{JSh69}.  In particular, if $R$ is regular and $x$ is in $\m^2$,
  then $S$ is not regular, so $\pd[S]{k}$ is infinite while $\pd{k}$
  is finite; see Theorem~A.
\end{example}

If $G$ is a totally reflexive $R$-module, then every $R$-regular
element is $G$-regular. A strong converse holds for modules of finite
projective dimension; it is (still) referred to as \emph{Auslander's
  zero-divisor conjecture:} %
\index{Auslander's zero-divisor conjecture} let $R$ be local and $M\ne
0$ be a finitely generated $R$-module with $\pd{M}$ finite. Then every
$M$-regular element is $R$-regular; for a proof see
Roberts~\cite[6.2.3]{PRb98}. An instance of \exaref{uv} shows that one
can \emph{not} relax the condition on $M$ to finite G-dimension:
 
\begin{example}
  \label{exa:4}
  Let $k$ be a field and consider the local ring $R =
  \pows[k]{X,Y}/(XY)$. Then the residue class $x$ of $X$ generates a
  totally reflexive module. The element $x$ is $(x)$-regular but
  nevertheless a zero-divisor in $R$.
\end{example}

While a tensor product of projective modules is projective, the next
example shows that totally reflexive modules do not have an analogous
property.

\begin{example}
  Let $R$ be as in \exaref{4}. The $R$-modules $(x)$ and $(y)$ are
  totally reflexive, but $\tp{(x)}{(y)} \is k$ is not. Indeed, $k$ is
  not a submodule of a free $R$-module.
\end{example}

The next result~\cite[4.13]{MAsMBr69} is parallel to Theorem~B in the
Introduction; it is known as the \emph{Auslander--Bridger Formula.} %
\index{Auslander--Bridger Formula}

\begin{theorem}
  \label{thm:ABr}
  Let $R$ be local. For every finitely generated $R$-module $M$ of
  finite G-dimension there is an equality $$\Gdim{M} = \dptR -
  \dpt{M}.$$
\end{theorem}

\begin{petit}
  \noindent
  In \cite{VMs00} Ma\c{s}ek corrects the proof of
  ~\cite[4.13]{MAsMBr69}. Proofs can also be found in
  \cite{MAs67}~and~\cite{lnm}.
\end{petit}

By \lemref{totref} the G-dimension is preserved under completion:

\begin{proposition}
  \label{prp:Gcpl}
  Let $R$ be local. For every finitely generated $R$-module $M$ there
  is an equality $$\Gdim{M} = \Gdim[\Rhat]{\tpp{\Rhat}{M}}.$$
\end{proposition}

The following main result from \cite[\S 3.2]{MAs67} is akin to
Theorem~A, but it differs in that it only deals with finitely
generated modules.

\begin{theorem}
  \label{thm:Gor-gdim}
  For a local ring $\Rmk$ the next conditions are equivalent.
  \begin{description}[$(iii)$]\setlength{\labelsep}{\mylabelsep}
  \item[\hfill $(i)$] $R$ is Gorenstein.%
    \index{Gorenstein ring!local}
  \item[\hfill $(ii)$] $\Gdim{k}$ is finite.
  \item[\hfill $(iii)$] $\Gdim{M}$ is finite for every finitely
    generated $R$-module $M$.
  \end{description}
\end{theorem}

\noindent 
It follows that non-projective totally reflexive modules exist over
any non-regular Gorenstein local ring. On the other hand, \exaref{taL}
shows that existence of such modules does not identify the ground ring
as a member of one of the standard classes, say, Cohen--Macaulay
rings.

A theorem of Christensen, Piepmeyer, Striuli, and
Takahashi~\cite[4.3]{CPST-08} shows that fewness of totally reflexive
modules comes in two distinct flavors:

\begin{theorem}
  Let $R$ be local. If there are only finitely many indecomposable
  totally reflexive $R$-modules, up to isomorphism, then $R$ is
  Gorenstein or every totally reflexive $R$-module is free.
\end{theorem}

\noindent This dichotomy brings two problems to light:

\begin{prb}
  \label{prb:CPST}
  Let $R$ be a local ring that is not Gorenstein and assume that there
  exists a non-free totally reflexive $R$-module. Find constructions
  that produce infinite families of non-isomorphic indecomposable
  totally reflexive modules.
\end{prb}

\begin{prb}
  \label{prb:LLAAMr}
  Describe the local rings over which every totally reflexive module
  is free.
\end{prb}

\noindent
While the first problem is posed in \cite{CPST-08}, the second one was
already raised by Avramov and Martsinkovsky~\cite{LLAAMr02}, and it is
proved \emph{ibid.}\ that over a Golod local ring that is not
Gorenstein, every totally reflexive module is free. Another partial
answer to \prbref{LLAAMr} is obtained by Yoshino~\cite{YYs03}, and by
Christensen and Veliche~\cite{LWCOVl07}.  The problem is also studied
by Takahashi in~\cite{RTk08a}.

Finally, Theorem~D in the Introduction motivates:

\begin{question}
  \label{qst:gdim-cm}
  Let $R$ be a local ring. If there exists a non-zero $R$-module of
  finite length and finite G-dimension, is then $R$ Cohen--Macaulay?%
  \index{Cohen--Macaulay local ring}
\end{question}

\noindent
A partial answer to this question is obtained by
Takahashi~\cite[2.3]{RTk04b}.

\subparagraph{Notes}

\begin{petit}
  A topic that was only treated briefly above is constructions of
  totally reflexive modules. Such constructions are found in
  \cite{AGP-97} by Avramov, Gasharov and Peeva, in work of Takahashi
  and Watanabe~\cite{RTkKWt07}, and in Yoshino's \cite{YYs03}.

  Hummel and Marley~\cite{LHmTMr} extend the notion of G-dimension to
  finitely presented modules over coherent rings and use it to define
  and study coherent Gorenstein~rings.

  Gerko~\cite[\S 2]{AAG01a} studies a dimension---the PCI-dimension or
  $\mathrm{CI}_*$-dimension---based on a subclass of the totally
  reflexive modules.  Golod~\cite{ESG84} studies a generalized notion
  of G-dimension: the $\mathrm{G}_C$-dimension, based on total
  reflexivity with respect to a semidualizing module~$C$. These
  studies are continued by, among others, Gerko~\cite[\S 1]{AAG01a},
  and Salarian, Sather-Wagstaff, and Yassemi~\cite{SSY-06}; see also
  the notes in \secref{fingdim}.

  An approach to homological dimensions that is not treated in this
  survey is based on so-called quasi-deformations. Several
  authors---among them Avramov, Gasharov, and Peeva \cite{AGP-97} and
  Veliche~\cite{OVl02}---take this approach to define homological
  dimensions that are intermediate between the projective dimension
  and the G-dimension for finitely generated modules. Gerko~\cite[\S
  3]{AAG01a} defines a Cohen--Macaulay dimension, which is a
  refinement of the G-dimension. Avramov \cite[\S 8]{LLA02} surveys
  these dimensions.
\end{petit}


\section{Gorenstein Projective Dimension}
\label{sec:gpd}

To extend the G-dimension beyond the realm of finitely generated
modules over Noetherian rings, Enochs and Jenda~\cite{EEnOJn95b}
introduced the notion of Gorenstein projective modules. The same
authors, and their collaborators, studied these modules in several
subsequent papers. The associated dimension, which is the focus of
this section, was studied by Christensen~\cite{lnm} and
Holm~\cite{HHl04a}.

In organization, this section is parallel to \secref{gdim}.

\begin{dfn}
  \label{dfn:GPm}
  An $R$-module $A$ is called \emph{Gorenstein projective} %
  \index{Gorenstein projective!module} if there exists an acyclic
  complex $\cx{P}$ of projective $R$-modules such that $\Coker{(P_1\to
    P_0)} \cong A$ and such that $\Hom{\cx{P}}{Q}$ is acyclic for
  every projective $R$-module $Q$.
\end{dfn}

\noindent
It is evident that every projective module is Gorenstein projective.

\begin{example}
  Every totally reflexive module is Gorenstein projective; this
  follows from \dfnref{taL} and \lemref{taL}.
\end{example}

Basic categorical properties are recorded in \cite[\S 2]{HHl04a}:

\begin{proposition}
  \label{prp:GPclosure}
  The class of Gorenstein projective $R$-modules is closed under
  direct sums and summands.
\end{proposition}

Every projective module is a direct summand of a free one. A parallel
result for Gorenstein projective modules, \thmref{SGPm} below, is due
to Bennis and Mahdou~\mbox{\cite[\S 2]{DBnNMh07}}; as substitute for
free modules they define:

\begin{dfn}
  \label{dfn:SGPm}
  An $R$-module $A$ is called \emph{strongly Gorenstein projective} if
  there exists an acyclic complex $\cx{P}$ of projective $R$-modules,
  in which all the differentials are identical, such that
  $\Coker{(P_1\to P_0)} \cong A$, and such that $\Hom{\cx{P}}{Q}$ is
  acyclic for every projective $R$-module $Q$.
\end{dfn}

\begin{theorem}
  \label{thm:SGPm}
  An $R$-module is Gorenstein projective if and only if it is a direct
  summand of a strongly Gorenstein projective $R$-module.
\end{theorem}

\begin{dfn}
  \label{dfn:GPr}
  An (augmented) \emph{Gorenstein projective resolution} of a module
  $M$ is an exact sequence \mbox{$\cdots \to A_i \to A_{i-1} \to
    \cdots \to A_0 \to M \to 0$}, where each module $A_i$ is
  Gorenstein projective.
\end{dfn}

\noindent
Note that every module has a Gorenstein projective resolution, as a
free resolution is trivially a Gorenstein projective one.

\begin{dfn}
  \label{dfn:Gpd}
  The \emph{Gorenstein projective dimension} %
  \index{Gorenstein projective!dimension} of a module \mbox{$M \ne
    0$}, denoted by $\Gpd{M}$, is the least integer \mbox{$n \geqslant
    0$} such that there exists a Gorenstein projective resolution of
  $M$ with \mbox{$A_i=0$} for all \mbox{$i>n$}.  If no such $n$
  exists, then $\Gpd{M}$ is infinite. By convention, set $\Gpd{0} =
  -\infty$.
\end{dfn}

In \cite[\S 2]{HHl04a} one finds the next standard theorem and
corollary.

\begin{theorem}
  \label{thm:Gpd}
  Let $M$ be an $R$-module of finite Gorenstein projective
  dimension. For every integer \mbox{$n \geqslant 0$} the following
  conditions are equivalent.
  \begin{description}[$(iii)$]\setlength{\labelsep}{\mylabelsep}
  \item[\hfill $(i)$] $\Gpd{M}\leqslant n$.
  \item[\hfill $(ii)$] $\Ext{i}{M}{Q}=0$ for all $i>n$ and all
    projective $R$-modules $Q$.
  \item[\hfill $(iii)$] $\Ext{i}{M}{N}=0$ for all $i>n$ and all
    $R$-modules $N$ with $\pd{N}$ finite.
  \item[\hfill $(iv)$] In every augmented Gorenstein projective
    resolution $$\cdots \to A_i \to A_{i-1} \to \cdots \to A_0 \to M
    \to 0$$ the module $\Coker{(A_{n+1} \to A_n)}$ is Gorenstein
    projective.
  \end{description}
\end{theorem}

\begin{corollary}
  \label{cor:measure-Gpd}
  For every $R$-module $M$ of finite Gorenstein projective dimension
  there is an equality $$\Gpd{M} = \sup\setof{\iinZ}{\Ext{i}{M}{Q} \ne
    0 \text{ for some projective $R$-module $Q$}}.$$
\end{corollary}

\begin{remark}
  For every $R$-module $M$ as in the corollary, the finitistic
  projective dimension of $R$ is an upper bound for $\Gpd{M}$;
  see~\cite[2.28]{HHl04a}.
\end{remark}

The next result \cite[2.24]{HHl04a} extends \prpref{3mg}.
   
\begin{proposition}
  \label{prp:3mgp}
  Let \mbox{$0 \to M' \to M \to M'' \to 0$} be an exact sequence of
  $R$-modules. If any two of the modules have finite Gorenstein
  projective dimension, then so has the~third.
\end{proposition}

The Gorenstein projective dimension is a refinement of the projective
dimension; this follows from \corref{measure-Gpd}:

\begin{proposition}
  \label{prp:Gpd-pd}
  For every $R$-module $M$ one has \mbox{$\Gpd{M} \leqslant \pd{M},$}
  and equality holds if $M$ has finite projective dimension.
\end{proposition}

\noindent
Supplementary information comes from Holm~\cite[2.2]{HHl04c}:

\begin{proposition}
  \label{prp:Gpd_equals_pd}
  If $M$ is an $R$-module of finite injective dimension, then there is
  an equality \mbox{$\Gpd{M}=\pd{M}$}.
\end{proposition}

The next result of Foxby is published in \cite[Ascent table
II(b)]{LWCHHl09a}.

\begin{proposition}
  \label{prp:GPbase}
  Let $S$ be an $R$-algebra of finite projective dimension. For every
  Gorenstein projective $R$-module $A$, the module $\tp{S}{A}$ is
  Gorenstein projective over $S$.
\end{proposition}

\subsubsection*{Noetherian Rings}
\label{sec:gpdn}

Finiteness of the Gorenstein projective dimension characterizes
Gorenstein rings. The next result of Enochs and
Jenda~\cite[12.3.1]{rha} extends \thmref{GorG}.

\begin{theorem}
  \label{thm:GP-Iwanaga-Gorenstein}
  \index{Gorenstein ring} Let $R$ be Noetherian and $n \ge 0$ be an
  integer. Then~$R$~is Gorenstein with $\dimR \le n$ if and only if
  \mbox{$\Gpd{M}\leqslant n$} for every $R$-module~$M$.
\end{theorem}

The next result \cite[4.2.6]{lnm} compares the Gorenstein projective
dimension to the G-dimension.

\begin{proposition}
  \label{prp:gdim-GPD}
  Let $R$ be Noetherian. For every finitely generated $R$-module $M$
  there is an equality $\Gpd{M}=\Gdim{M}$.
\end{proposition}

The Gorenstein projective dimension of a module can not be measured
locally; that is, \prpref{gdimp} does not extend to non-finitely
generated modules.  As a consequence of \prpref{GPbase}, though, one
has the following:

\begin{proposition}
  Let $R$ be Noetherian of finite Krull dimension. For every
  $R$-module $M$ and every prime ideal $\p$ in $R$ one has
  $\Gpd[R_\p]{M_\p} \leqslant \Gpd{M}$.
\end{proposition}

Theorem~E and \prpref{Gpd_equals_pd} yield:

\begin{theorem}
  \label{thm:Gorenstein_on_Supp1}
  Let $R$ be Noetherian and $M$ a finitely generated $R$-module. If
  $\Gpd{M}$ and $\id{M}$ are finite, then $R_\p$ is Gorenstein for all
  $\p \in \Supp{M}$.
\end{theorem}

\subsubsection*{Local Rings}
\label{sec:gpdl}

The next characterization of Gorenstein local rings---akin to
Theorem~A in the Introduction---follows from
\thmref[Theorems~]{Gor-gdim} and \thmref[]{GP-Iwanaga-Gorenstein} via
\prpref{gdim-GPD}.

\begin{theorem}
  \label{thm:Gor-Gpd}
  For a local ring $\Rmk$ the next conditions are equivalent.
  \begin{description}[$(iii)$]\setlength{\labelsep}{\mylabelsep}
  \item[\hfill $(i)$] $R$ is Gorenstein.%
    \index{Gorenstein ring!local}
  \item[\hfill $(ii)$] $\Gpd{k}$ is finite.
  \item[\hfill $(iii)$] $\Gpd{M}$ is finite for every $R$-module $M$.
  \end{description}
\end{theorem}

The inequality in the next theorem is a consequence of
\prpref{GPbase}. The second assertion is due to Esmkhani and
Tousi~\cite[3.5]{MAEMTs07}, cf.~\cite[4.1]{CFH-06}. The result should
be compared to \prpref{Gcpl}.

\begin{theorem}
  \label{thm:GPcompletion}
  Let $R$ be local and $M$ be an $R$-module. Then one
  has $$\Gpd[\Rhat]{(\tp{\Rhat}{M})} \le \Gpd{M},$$ and if
  $\Gpd[\Rhat]{(\tp{\Rhat}{M})}$ is finite, then so is $\Gpd{M}$.
\end{theorem}

\subparagraph{Notes}

\begin{petit}
  Holm and J\o rgensen~\cite{HHlPJr06} extend Golod's \cite{ESG84}
  notion of $\mathrm{G}_C$-dimension to non-finitely generated modules
  in the form of a $C$-Gorenstein projective dimension. Further
  studies of this dimension are made by White~\cite{DWh}.
\end{petit}

\section{Gorenstein Injective Dimension}
\label{sec:gid}

The notion of Gorenstein injective modules is (categorically) dual to
that of Gorenstein projective modules. The two were introduced in the
same paper by Enochs and Jenda~\cite{EEnOJn95b} and investigated in
subsequent works by the same authors, by Christensen and
Sather-Wagstaff~\cite{LWCSSWc}, and by Holm~\cite{HHl04a}.

This section is structured parallelly to the previous ones.

\begin{dfn}
  \label{dfn:GIm}
  An $R$-module $B$ is called \emph{Gorenstein injective} %
  \index{Gorenstein injective!module} if there exists an acyclic
  complex $\cx{I}$ of injective $R$-modules such that $\Ker{(I^0\to
    I^1)} \cong B$, and such that $\Hom{E}{\cx{I}}$ is acyclic for
  every injective $R$-module $E$.
\end{dfn}

\noindent
It is clear that every injective module is Gorenstein injective.

\begin{example}
  \label{exa:gim}
  Let $\cx{L}$ be a totally acyclic complex of finitely generated
  projective $R$-modules, see~\dfnref{taL}, and let $I$ be an
  injective $R$-module. Then the acyclic complex $\cx{I} =
  \Hom{\cx{L}}{I}$ consists of injective modules, and it follows from
  \lemref{taL} that the complex $\Hom{E}{\cx{I}} \is
  \Hom{\tp{E}{\cx{L}}}{I}$ is acyclic for every injective module
  $E$. Thus, if $G$ is a totally reflexive $R$-module and $I$ is
  injective, then the module $\Hom{G}{I}$ is Gorenstein injective.
\end{example}

Basic categorical properties are established in \cite[2.6]{HHl04a}:

\begin{proposition}
  \label{prp:GIclosure}
  The class of Gorenstein injective $R$-modules is closed under direct
  products and summands.
\end{proposition}

\noindent
Under extra assumptions on the ring, \thmref{GIcolimit} gives more
information.

\begin{dfn}
  \label{dfn:GIr}
  An (augmented) \emph{Gorenstein injective resolution} of a module
  $M$ is an exact sequence \mbox{$0 \to M \to B^0 \to \cdots \to
    B^{i-1} \to B^i \to \cdots$}, where each module $B^i$ is
  Gorenstein injective.
\end{dfn}

\noindent
Note that every module has a Gorenstein injective resolution, as an
injective resolution is trivially a Gorenstein injective one.

\begin{dfn}
  \label{dfn:Gid}
  The \emph{Gorenstein injective dimension} %
  \index{Gorenstein injective!dimension} of an $R$-module \mbox{$M \ne
    0$}, denoted by $\Gid{M}$, is the least integer \mbox{$n \geqslant
    0$} such that there exists a Gorenstein injective resolution of
  $M$ with \mbox{$B^i=0$} for all \mbox{$i>n$}.  If no such $n$
  exists, then $\Gid{M}$ is infinite. By convention, set $\Gid{0} =
  -\infty$.
\end{dfn}

The next standard theorem is \cite[2.22]{HHl04a}.

\begin{theorem}
  \label{thm:Gid}
  Let $M$ be an $R$-module of finite Gorenstein injective
  dimension. For every integer \mbox{$n \geqslant 0$} the following
  conditions are equivalent.
  \begin{description}[$(iii)$]\setlength{\labelsep}{\mylabelsep}
  \item[\hfill $(i)$] $\Gid{M}\leqslant n$.
  \item[\hfill $(ii)$] $\Ext{i}{E}{M}=0$ for all $i>n$ and all
    injective $R$-modules $E$.
  \item[\hfill $(iii)$] $\Ext{i}{N}{M}=0$ for all $i>n$ and all
    $R$-modules $N$ with $\id{N}$ finite.
  \item[\hfill $(iv)$] In every augmented Gorenstein injective
    resolution $$0 \to M \to B^0 \to \cdots \to B^{i-1} \to B^i \to
    \cdots$$ the module $\Ker{(B^n \to B^{n+1})}$ is Gorenstein
    injective.
  \end{description}
\end{theorem}

\begin{corollary}
  \label{cor:measure-Gid}
  For every $R$-module $M$ of finite Gorenstein injective dimension
  there is an equality $$\Gid{M} = \sup\setof{\iinZ}{\Ext{i}{E}{M} \ne
    0 \text{ for some injective $R$-module $E$}}.$$
\end{corollary}

\begin{remark}
  For every $R$-module $M$ as in the corollary, the finitistic
  injective dimension of $R$ is an upper bound for $\Gid{M}$;
  see~\cite[2.29]{HHl04a}.
\end{remark}

The next result \cite[2.25]{HHl04a} is dual to \prpref{3mgp}.

\begin{proposition}
  \label{prp:3Ginj}
  Let \mbox{$0 \to M' \to M \to M'' \to 0$} be an exact sequence of
  $R$-modules. If any two of the modules have finite Gorenstein
  injective dimension, then so has the~third.
\end{proposition}

The Gorenstein injective dimension is a refinement of the injective
dimension; this follows from \corref{measure-Gid}:

\begin{proposition}
  \label{prp:Gid-id}
  For every $R$-module $M$ one has \mbox{$\Gid{M} \leqslant \id{M}$},
  and equality holds if $M$ has finite injective dimension.
\end{proposition}

\noindent
Supplementary information comes from Holm~\cite[2.1]{HHl04c}:

\begin{proposition}
  \label{prp:Gid_equals_id}
  If $M$ is an $R$-module of finite projective dimension, then there
  is an equality \mbox{$\Gid{M}=\id{M}$}. In particular, one has
  $\Gid{R}=\id{R}$.
\end{proposition}

In \cite{LWCHHl09a} Christensen and Holm study (co)base change of
modules of finite Gorenstein homological dimension. The following is
elementary to verify:

\begin{proposition}
  \label{prp:GIcobase}
  Let $S$ be an $R$-algebra of finite projective dimension. For every
  Gorenstein injective $R$-module $B$, the module $\Hom{S}{B}$ is
  Gorenstein injective over $S$.
\end{proposition}

\noindent
For a conditional converse see \thmref[Theorems~]{GI_cotorsion} and
\thmref[]{GidCoBase}.

The next result of Bennis, Mahdou, and Ouarghi~\cite[2.2]{BMO} should
be compared to characterizations of Gorenstein rings like
\thmref[Theorems~]{GP-Iwanaga-Gorenstein} and
\thmref[]{Iwanaga-Gorenstein}, and also to
\thmref[Theorems~]{Gorenstein_on_Supp1} and
\thmref[]{Gorenstein_on_Supp}. It is a perfect Gorenstein counterpart
to a classical result due to Faith and Walker among others; see
e.g.~\cite[4.2.4]{wei}.

\begin{theorem}
  \label{thm:qfi}
  The following conditions on $R$ are equivalent.
  \begin{description}[$(iii)$]\setlength{\labelsep}{\mylabelsep}
  \item[\hfill $(i)$] $R$ is quasi-Frobenius. %
    \index{quasi-Frobenius ring}
  \item[\hfill $(ii)$] Every $R$-module is Gorenstein projective.
  \item[\hfill $(iii)$] Every $R$-module is Gorenstein injective.
  \item[\hfill $(iv)$] Every Gorenstein projective $R$-module is
    Gorenstein injective.
  \item[\hfill $(v)$] Every Gorenstein injective $R$-module is
    Gorenstein projective.
  \end{description}
\end{theorem}

\subsubsection*{Noetherian Rings}
\label{sec:gidn}

Finiteness of the Gorenstein injective dimension characterizes
Gorenstein rings; this result is due to Enochs and
Jenda~\cite[3.1]{EEnOJn94b}:

\begin{theorem}
  \label{thm:Iwanaga-Gorenstein}
  Let $R$ be Noetherian and $n \ge 0$ be an integer. Then $R$ is
  Gorenstein with $\dimR \le n$ if and only if
  \mbox{$\,\Gid{M}\leqslant n$} for every $R$-module~$M$.%
  \index{Gorenstein ring}
\end{theorem}

A ring is Noetherian if every countable direct sum of injective
modules is injective (and only if every direct limit of injective
modules is injective). The ``if'' part has a perfect Gorenstein
counterpart:

\begin{proposition}
  \label{prp:count}
  If every countable direct sum of Gorenstein injective $R$-modules is
  Gorenstein injective, then $R$ is Noetherian.
\end{proposition}

\begin{petit}
  \begin{proof}
    It is sufficient to see that every countable direct sum of
    injective $R$-modules is injective. Let $\{E_n\}_{n \in \NN}$ be a
    family of injective modules. By assumption, the module $\bigoplus
    E_n$ is Gorenstein injective; in particular, there is an
    epimorphism $\mapdef[\twoheadrightarrow]{\f}{I}{\bigoplus E_n}$
    such that $I$ is injective and $\Hom{E}{\f}$ is surjective for
    every injective $R$-module $E$. Applying this to $E=E_n$ it is
    elementary to verify that $\f$ is a split epimorphism. \qed
  \end{proof}
\end{petit}

Christensen, Frankild, and Holm~\cite[6.9]{CFH-06} provide a partial
converse:

\begin{theorem}
  \label{thm:GIcolimit}
  Assume that $R$ is a homomorphic image of a Gorenstein ring of
  finite Krull dimension. Then the class of Gorenstein injective
  modules is closed under direct limits; in particular, it is closed
  under direct sums.
\end{theorem}

\noindent
As explained in the Introduction, the hypothesis on $R$ in this
theorem ensures the existence of a dualizing $R$-complex and an
associated Bass class, cf.~\secref{fingd}. These tools are essential
to the known proof of \thmref{GIcolimit}.

\begin{question}
  Let $R$ be Noetherian. Is then every direct limit of Gorenstein
  injective $R$-modules Gorenstein injective?
\end{question}

Next follows a Gorenstein version of Chouinard's formula
\cite[3.1]{LGC76}; %
\index{Chouinard Formula!for Gid} it is proved
in~\cite[2.2]{LWCSSWc}. Recall that the \emph{width} %
\index{width} of a module $M$ over a local ring $\Rmk$ is defined
as $$\wdt{M} = \inf\setof{i\in\ZZ}{\Tor{i}{k}{M}\ne 0}.$$

\begin{theorem}
  \label{thm:Gid-Chouinard}
  Let $R$ be Noetherian. For every $R$-module $M$ of finite Gorenstein
  injective dimension there is an equality
  \begin{displaymath}
    \Gid{M}=\setof{\dpt[]{R_\p}-\wdt[R_\p]{M_\p}}{\p \in \SpecR}.
  \end{displaymath}
\end{theorem}

Let $M$ be an $R$-module, and let $\p$ be a prime ideal in
$R$. Provided that $\Gid[R_\p]{M_\p}$ is finite, the inequality
$\Gid[R_\p]{M_\p} \le \Gid{M}$ follows immediately from the
theorem. However, the next question remains open.

\begin{question}
  Let $R$ be Noetherian and $B$ be a Gorenstein injective
  $R$-mod\-ule. Is then $B_\p$ Gorenstein injective over $R_\p$ for
  every prime ideal $\p$ in $R\,$?
\end{question}

\noindent
A partial answer is known from \cite[5.5]{CFH-06}:

\begin{proposition}
  Assume that $R$ is a homomorphic image of a Gorenstein ring of
  finite Krull dimension. For every $R$-module $M$ and every prime
  ideal~$\p$ there is an inequality
  \mbox{$\Gid[R_\p]{M_\p}\leqslant\Gid{M}.$}
\end{proposition}

Theorem~E and \prpref{Gid_equals_id} yield:

\begin{theorem}
  \label{thm:Gorenstein_on_Supp}
  Let $R$ be Noetherian and $M$ a finitely generated $R$-mod\-ule. If
  $\Gid{M}$ and $\pd{M}$ are finite, then $R_\p$ is Gorenstein for all
  $\p \in \Supp{M}$.
\end{theorem}

\subsubsection*{Local Rings}
\label{sec:gidl}

The following theorem of Foxby and Frankild~\cite[4.5]{HBFAJF07}
generalizes work of Peskine and Szpiro \cite{CPsLSz73}, cf.~Theorem~E.

\begin{theorem}
  \label{thm:localGorenstein}
  A local ring $R$ is Gorenstein if and only if there exists a
  non-zero cyclic $R$-module of finite Gorenstein injective
  dimension.%
  \index{Gorenstein ring!local}
\end{theorem}

\thmref[Theorems~]{Iwanaga-Gorenstein} and \thmref[]{localGorenstein}
yield a parallel to \thmref[Thm.~]{Gor-gdim}, akin to Theorem~A.

\begin{corollary}
  \label{cor:Gor-Gid}
  For a local ring $\Rmk$ the next conditions are equivalent.
  \begin{description}[$(iii)$]\setlength{\labelsep}{\mylabelsep}
  \item[\hfill $(i)$] $R$ is Gorenstein.%
    \index{Gorenstein ring!local}
  \item[\hfill $(ii)$] $\Gid{k}$ is finite.
  \item[\hfill $(iii)$] $\Gid{M}$ is finite for every $R$-module $M$.
  \end{description}
\end{corollary}

The first part of the next theorem is due to Christensen, Frankild,
and Iyengar, and published in~\cite[3.6]{HBFAJF07}. The equality in
\thmref{Bass_formula}---the Gorenstein analogue of Theorem~C in the
Introduction---is proved by Khatami, Tousi, and
Yassemi~\cite[2.5]{KTY-}; see also \cite[2.3]{LWCSSWc}.

\begin{theorem}
  \label{thm:Bass_formula} %
  \index{Bass Formula!for Gid} Let $R$ be local and \mbox{$M\ne 0$} be
  a finitely generated $R$-module. Then $\Gid{M}$ and
  $\Gid[\Rhat]{\tpp{\Rhat}{M}}$ are simultaneously finite, and when
  they are finite, there is an equality $$\Gid{M}=\dptR.$$
\end{theorem}

\begin{remark}
  Let $R$ be local and $M\ne 0$ be an $R$-module. If $M$ has finite
  length and finite G-dimension, then its Matlis dual has finite
  Gorenstein injective dimension, cf.~\exaref{gim}. See also
  Takahashi~\cite{RTk06a}.
\end{remark}

\noindent
This remark and Theorem~D from the Introduction motivate:

\begin{question}
  \label{qst:gid-cm}
  Let $R$ be local. If there exists a non-zero finitely generated
  $R$-module of finite Gorenstein injective dimension, is then $R$
  Cohen--Macaulay?%
  \index{Cohen--Macaulay local ring}
\end{question}

\noindent
A partial answer to this question is given by
Yassemi~\cite[1.3]{SYs07}.

Esmkhani and Tousi~\cite[2.5]{MAEMTs07a} prove the following
conditional converse to \prpref{GIcobase}. Recall that an $R$-module
$M$ is said to be \emph{cotorsion} %
\index{cotorsion module} if $\Ext{1}{F}{M}=0$ for every flat
$R$-module $F$.

\begin{theorem}
  \label{thm:GI_cotorsion}
  Let $R$ be local. An $R$-module $M$ is Gorenstein injective if and
  only if it is cotorsion and $\Hom{\Rhat}{M}$ is Gorenstein injective
  over $\Rhat$.
\end{theorem}

\noindent
The example below demonstrates the necessity of the cotorsion
hypothesis. Working in the derived category one obtains a stronger
result; see~\thmref[Thm.~]{GidCoBase}.

\begin{example}
  \label{exa:GIcobase}
  Let $\Rm$ be a local domain which is not $\m$-adically
  complete. Aldrich, Enochs, and L\'opez-Ramos~\cite[3.3]{AEL-02} show
  that the module $\Hom{\Rhat}{R}$ is zero and hence Gorenstein
  injective over $\Rhat$. However, $\Gid{R}$ is infinite if $R$ is not
  Gorenstein, cf.~\prpref{Gid_equals_id}.
\end{example}

\subparagraph{Notes}

\begin{petit}
  Dual to the notion of strongly Gorenstein projective modules, see
  \dfnref{SGPm}, Bennis and Mahdou~\cite{DBnNMh07} also study strongly
  Gorenstein injective modules.

  Several authors---Asadollahi, Sahandi, Salarian, Sazeedeh, Sharif,
  and Yassemi---have studied the relationship between Gorenstein
  injectivity and local cohomology; see \cite{JAsSSl06},
  \cite{PShTSH07}, \cite{RSz04}, \cite{RSz07}, and \cite{SYs07}.
\end{petit}


\section{Gorenstein Flat Dimension}
\label{sec:gfd}

Another extension of the G-dimension is based on Gorenstein flat
modules---a notion due to Enochs, Jenda, and
Torrecillas~\cite{EJT-93}.  Christensen \cite{lnm} and
Holm~\cite{HHl04a} are other main references for this section.

The organization of this section follows the pattern from
\secref[Sections~]{gdim}--\secref[]{gid}.

\begin{dfn}
  \label{dfn:GFm}
  \index{Gorenstein flat!module} An $R$-module $A$ is called
  \emph{Gorenstein flat} if there exists an acyclic complex $\cx{F}$
  of flat $R$-modules such that $\Coker{(F_1\to F_0)} \cong A$, and
  such that $\tp{E}{\cx{F}}$ is acyclic for every injective $R$-module
  $E$.
\end{dfn}

\noindent
It is evident that every flat module is Gorenstein flat.

\begin{example}
  \label{exa:ta-implies-GF}
  Every totally reflexive module is Gorenstein flat; this follows from
  \dfnref{taL} and \lemref{taL}.
\end{example}

Here is a direct consequence of \dfnref{GFm}:

\begin{proposition}
  \label{prp:GFclosure}
  The class of Gorenstein flat $R$-modules is closed under direct
  sums.
\end{proposition}

\noindent
See \thmref[Theorems~]{GFprod} and \thmref[]{GFclosureII} for further
categorical properties of Gorenstein flat modules.

\begin{dfn}
  \label{dfn:GFr}
  An (augmented) \emph{Gorenstein flat resolution} of a module $M$ is
  an exact sequence \mbox{$\cdots \to A_i \to A_{i-1} \to \cdots \to
    A_0 \to M \to 0$}, where each module $A_i$ is Gorenstein flat.
\end{dfn}

\noindent
Note that every module has a Gorenstein flat resolution, as a free
resolution is trivially a Gorenstein flat one.

\begin{dfn}
  \label{dfn:Gfd}
  The \emph{Gorenstein flat dimension} %
  \index{Gorenstein flat!dimension} of an $R$-module \mbox{$M \ne 0$},
  denoted by $\Gfd{M}$, is the least integer \mbox{$n \geqslant 0$}
  such that there exists a Gorenstein flat resolution of $M$ with
  \mbox{$A_i=0$} for all \mbox{$i>n$}.  If no such $n$ exists, then
  $\Gfd{M}$ is infinite. By convention, set $\Gfd{0} = -\infty$.
\end{dfn}

The next duality result is an immediate consequence of the
definitions.

\begin{proposition}
  \label{prp:GFdual}
  Let $M$ be an $R$-module. For every injective $R$-module $E$ there
  is an inequality $\Gid{\Hom{M}{E}} \leqslant \Gfd{M}$.
\end{proposition}

\noindent

Recall that an $R$-module $E$ is called \emph{faithfully} injective if
it is injective and \mbox{$\Hom{M}{E}=0$} only if \mbox{$M=0$}. The
next question is inspired by the classical situation. It has an
affirmative answer for Noetherian rings; see \thmref{GFdual}.%
\index{faithfully injective module}

\begin{question}
  \label{qst:HomME}
  Let $M$ and $E$ be $R$-modules. If $E$ is faithfully injective and
  the module $\Hom{M}{E}$ is Gorenstein injective, is then $M$
  Gorenstein flat?
\end{question}

A straightforward application of \prpref{GFdual} shows that the
Gorenstein flat dimension is a refinement of the flat dimension; cf.\
Bennis~\cite[2.2]{DBn}:

\begin{proposition}
  \label{prp:Gfd-fd}
  For every $R$-module $M$ one has \mbox{$\Gfd{M} \leqslant \fd{M}$},
  and equality holds if $M$ has finite flat dimension.
\end{proposition}

The following result is an immediate consequence of \dfnref{GFm}. Over
a local ring a stronger result is available; see
\thmref{GFcompletion}.

\begin{proposition}
  \label{prp:GFbase}
  Let $S$ be a flat $R$-algebra. For every $R$-module $M$ there is an
  inequality $\Gfd[S]{(\tp{S}{M})} \leqslant \Gfd{M}$.
\end{proposition}

\begin{corollary}
  Let $M$ be an $R$-module. For every prime ideal $\p$ in $R$ there is
  an inequality $\Gfd[R_\p]{M_\p} \leqslant \Gfd{M}$.
\end{corollary}

\subsubsection*{Noetherian Rings}
\label{sec:gfdn}

Finiteness of the Gorenstein flat dimension characterizes Gorenstein
rings; this is a result of Enochs and Jenda~\cite[3.1]{EEnOJn94b}:

\begin{theorem}
  \label{thm:GF-Iwanaga-Gorenstein}
  Let $R$ be Noetherian and $n \ge 0$ be an integer. Then $R$ is
  Gorenstein with $\dimR \le n$ if and only if
  \mbox{$\,\Gfd{M}\leqslant n$} for every $R$-module~$M$.%
  \index{Gorenstein ring}
\end{theorem}

A ring is coherent if and only if every direct product of flat modules
is flat. We suggest the following problem:

\begin{prb}
  Describe the rings over which every direct product of Gorenstein
  flat modules is Gorenstein flat.
\end{prb}

\noindent
Partial answers are due to Christensen, Frankild, and
Holm~\cite[5.7]{CFH-06} and to Murfet and Salarian
\cite[6.21]{DMfSSl}.

\begin{theorem}
  \label{thm:GFprod}
  Let $R$ be Noetherian. The class of Gorenstein flat $R$-modules is
  closed under direct products under either of the following
  conditions:
  \begin{description}[$(a)$]\setlength{\labelsep}{\mylabelsep}
  \item[\hfill\rm (a)] $R$ is homomorphic image of a Gorenstein ring
    of finite Krull dimension.
  \item[\hfill\rm (b)] $R_\p$ is Gorenstein for every non-maximal
    prime ideal $\p$ in $R$.
  \end{description}
\end{theorem}

The next result follows from work of Enochs, Jenda, and L\'opez-Ramos
\cite[2.1]{EJL-04} and \cite[3.3]{EEnOJn85}.

\begin{theorem}
  \label{thm:GFclosureII}
  Let $R$ be Noetherian. Then the class of Gorenstein flat $R$-modules
  is closed under direct summands and direct limits.
\end{theorem}

A result of Govorov~\cite{VEG65} and Lazard \cite[1.2]{DLz69} asserts
that a module is flat if and only if it is a direct limit of finitely
generated projective modules. For Gorenstein flat modules, the
situation is more complicated:

\begin{remark}
  \label{rmk:Lazard}
  Let $R$ be Noetherian. It follows from \exaref{ta-implies-GF} and
  \thmref{GFclosureII} that a direct limit of totally reflexive
  modules is Gorenstein flat. If $R$ is Gorenstein of finite Krull
  dimension, then every Gorenstein flat $R$-module can be written as a
  direct limit of totally reflexive modules; see Enochs and
  Jenda~\cite[10.3.8]{rha}. If $R$ is not Gorenstein, this conclusion
  may fail; see Beligiannis and Krause~\cite[4.2 and 4.3]{ABgHKr} and
  \thmref{HHlPJr}.
\end{remark}

The next result \cite[6.4.2]{lnm} gives a partial answer to
\qstref{HomME}.

\begin{theorem}
  \label{thm:GFdual}
  Let $R$ be Noetherian, and let $M$ and $E$ be $R$-modules. If $E$ is
  faithfully injective, then there is an equality $$\Gid{\Hom{M}{E}} =
  \Gfd{M}.$$
\end{theorem}

\thmref{Gfd} is found in \cite[3.14]{HHl04a}. It can be obtained by
application of \thmref{GFdual} to \thmref{Gid}.

\begin{theorem}
  \label{thm:Gfd}
  Let $R$ be Noetherian and $M$ be an $R$-module of finite Gorenstein
  flat dimension. For every integer \mbox{$n \geqslant 0$} the
  following are~equivalent.
  \begin{description}[$(iii)$]\setlength{\labelsep}{\mylabelsep}
  \item[\hfill $(i)$] $\Gfd{M}\leqslant n$.
  \item[\hfill $(ii)$] $\Tor{i}{E}{M}=0$ for all $i>n$ and all
    injective $R$-modules $E$.
  \item[\hfill $(iii)$] $\Tor{i}{N}{M}=0$ for all $i>n$ and all
    $R$-modules $N$ with $\id{N}$ finite.
  \item[\hfill $(iv)$] In every augmented Gorenstein flat resolution
    $$\cdots \to A_i \to A_{i-1} \to \cdots \to A_0 \to M \to 0$$
    the module $\Coker{(A_{n+1} \to A_n)}$ is Gorenstein flat.
  \end{description}
\end{theorem}

\begin{corollary}
  \label{cor:measure-Gfd}
  Let $R$ be Noetherian.  For every $R$-module $M$ of finite
  Gorenstein flat dimension there is an equality $$\Gfd{M} =
  \sup\setof{\iinZ}{\Tor{i}{E}{M} \ne 0 \text{ for some injective
      $R$-module $E$}}.$$
\end{corollary}

\begin{remark}
  For every $R$-module $M$ as in the corollary, the finitistic flat
  dimension of $R$ is an upper bound for $\Gfd{M}$;
  see~\cite[3.24]{HHl04a}.
\end{remark}

The next result \cite[3.15]{HHl04a} follows by \thmref{GFdual} and
\prpref{3Ginj}.

\begin{proposition}
  \label{prp:3gf}
  Let $R$ be Noetherian. If any two of the modules in an exact
  sequence \mbox{$0 \to M' \to M \to M'' \to 0$} have finite
  Gorenstein flat dimension, then so has the third.
\end{proposition}

A result of Holm~\cite[2.6]{HHl04c} supplements \prpref{Gfd-fd}:

\begin{proposition}
  \label{prp:Gfd_equals_fd}
  Let $R$ be Noetherian of finite Krull dimension.  For every
  $R$-module $M$ of finite injective dimension one has
  \mbox{$\Gfd{M}=\fd{M}$}.
\end{proposition}

Recall that the \emph{depth} %
\index{depth} of a module $M$ over a local ring $\Rmk$ is given
as $$\dpt{M} = \inf\setof{i\in\ZZ}{\Ext{i}{k}{M}\ne 0}.$$
\thmref{Gfd-Chouinard} is a Gorenstein version of Chouinard's
\cite[1.2]{LGC76}. %
\index{Chouinard Formula!for Gfd} It follows from \cite[3.19]{HHl04a}
and \cite[2.4(b)]{CFF-02}; see also Iyengar and
Sather-Wagstaff~\cite[8.8]{SInSSW04}.

\begin{theorem}
  \label{thm:Gfd-Chouinard}
  Let $R$ be Noetherian. For every $R$-module $M$ of finite Gorenstein
  flat dimension there is an equality
  \begin{displaymath}
    \Gfd{M} = \setof{\dpt[]{R_\p}-\dpt[R_\p]{M_\p}}{\p \in \SpecR}.
  \end{displaymath}
\end{theorem}

The next two results compare the Gorenstein flat dimension to the
Gorenstein projective dimension. The inequality in \thmref{Gfd-Gpd} is
\cite[3.4]{HHl04a}, and the second assertion in this theorem is due to
Esmkhani and Tousi~\cite[3.4]{MAEMTs07}.
\pagebreak[4]

\begin{theorem}
  \label{thm:Gfd-Gpd}
  Let $R$ be Noetherian of finite Krull dimension, and let $M$ be an
  $R$-module.  Then there is an inequality $$\Gfd{M} \leqslant
  \Gpd{M},$$ and if $\Gfd{M}$ is finite, then so is $\Gpd{M}$.
\end{theorem}

\noindent
It is not known whether the inequality in \thmref{Gfd-Gpd} holds over
every commutative ring. For finitely generated modules one has
\cite[4.2.6 and 5.1.11]{lnm}:

\begin{proposition}
  \label{prp:Gfd-gdim}
  Let $R$ be Noetherian. For every finitely generated $R$-module $M$
  there is an equality $\Gfd{M}=\Gpd{M}=\Gdim{M}$.
\end{proposition}

The next result \cite[5.1]{CFH-06} is related to \thmref{GFdual}; the
question that follows is prompted by the classical situation.

\begin{theorem}
  \label{thm:noname}
  Assume that $R$ is a homomorphic image of a Gorenstein ring of
  finite Krull dimension. For every $R$-module $M$ and every injective
  $R$-module $E$ there is an inequality $$\Gfd{\Hom{M}{E}} \leqslant
  \Gid{M},$$ and equality holds if $E$ is faithfully injective.
\end{theorem}

\begin{question}
  Let $R$ be Noetherian and $M$ and $E$ be $R$-modules. If $M$ is
  Gorenstein injective and $E$ is injective, is then $\Hom{M}{E}$
  Gorenstein flat?
\end{question}

\subsubsection*{Local Rings}
\label{sec:gfdl}

Over a local ring there is a stronger version~\cite[3.5]{MAEMTs07} of
\prpref{GFbase}:

\begin{theorem}
  \label{thm:GFcompletion}
  Let $R$ be local. For every $R$-module $M$ there is an equality
  $$\Gfd[\Rhat](\tp{\Rhat}{M})=\Gfd{M}.$$
\end{theorem}

Combination of~\cite[2.1 and 2.2]{HHl04c} with Theorem~E yields the
next result. Recall that a non-zero finitely generated module has
finite depth.

\begin{theorem}
  \label{thm:Gor-Gfd}
  For a local ring $R$ the following conditions are equivalent.
  \begin{description}[$(iii)$]\setlength{\labelsep}{\mylabelsep}
  \item[\hfill $(i)$] $R$ is Gorenstein.%
    \index{Gorenstein ring!local}
  \item[\hfill $(ii)$] There is an $R$-module $M$ with $\dpt{M}$,
    $\fd{M}$, and $\Gid{M}$ finite.
  \item[\hfill $(iii)$] There is an $R$-module $M$ with $\dpt{M}$,
    $\id{M}$, and $\Gfd{M}$ finite.
  \end{description}
\end{theorem}

\noindent
We have for a while been interested in:

\begin{question}
  \label{qst:GHBF}
  Let $R$ be local. If there exists an $R$-module $M$ with $\dpt{M}$,
  $\Gfd{M}$, and $\Gid{M}$ finite, is then $R$ Gorenstein?
\end{question}

A theorem of J{\o}rgensen and Holm~\cite{HHlPJr} brings perspective to
\rmkref{Lazard}.

\begin{theorem}
  \label{thm:HHlPJr}
  Assume that $R$ is Henselian local and a homomorphic image of a
  Gorenstein ring. If every Gorenstein flat $R$-module is a direct
  limit of totally reflexive modules, then $R$ is Gorenstein or every
  totally reflexive $R$-module is free.
\end{theorem}

\subparagraph{Notes}

\begin{petit}
  Parallel to the notion of strongly Gorenstein projective modules,
  see \dfnref{SGPm}, Bennis and Mahdou~\cite{DBnNMh07} also study
  strongly Gorenstein flat modules. A different notion of strongly
  Gorenstein flat modules is studied by Ding, Li, and Mao in
  \cite{DLM-}.
\end{petit}


\section{Relative Homological Algebra}
\label{sec:rha}

Over a Gorenstein local ring, the totally reflexive modules are
exactly the maximal Cohen--Macaulay modules, and their representation
theory is a classical topic. Over rings that are not Gorenstein, the
representation theory of totally reflexive modules was taken up by
Takahashi~\cite{RTk05a} and Yoshino~\cite{YYs03}. Conclusive results
have recently been obtained by Christensen, Piepmeyer, Striuli, and
Takahashi~\cite{CPST-08} and by Holm and J{\o}rgensen~\cite{HHlPJr}.
These results are cast in the language of precovers and preenvelopes;
see \thmref{Gpre}.

Relative homological algebra studies dimensions and (co)homology
functors based on resolutions that are constructed via precovers or
preenvelopes. Enochs and Jenda and their collaborators have made
extensive studies of the precovering and preenveloping properties of
the classes of Gorenstein flat and Gorenstein injective modules. Many
of their results are collected in~\cite{rha}.

\subsubsection*{Terminology}

Let $\cH$ be a class of $R$-modules. Recall that an
\emph{$\cH$-precover} %
\index{precover} (also called a right $\cH$-approximation) %
\index{approximation} of an $R$-module $M$ is a homomorphism
$\mapdef{\f}{H}{M}$ with $H$ in $\cH$ such
that $$\dmapdef{\Hom{H'}{\f}}{\Hom{H'}{H}}{\Hom{H'}{M}}$$ is
surjective for every $H'$ in $\cH$. That is, every homomorphism from a
module in $\cH$ to $M$ factors through $\f$. Dually one defines
$\cH$-preenvelopes (also called left $\cH$-approximations).

\begin{remark}
  If $\cH$ contains all projective modules, then every $\cH$-precover
  is an epimorphism. Thus, Gorenstein projective/flat precovers are
  epimorphisms.

  If $\cH$ contains all injective modules, then every
  $\cH$-preenvelope is a monomorphism. Thus, every Gorenstein
  injective preenvelope is a monomorphism.
\end{remark}

Fix an $\cH$-precover $\f$. It is called \emph{special} if the
cohomology module $\Ext{1}{H'}{\Ker{\f}}$ vanishes for every module
$H'$ in $\cH$. It is called a \emph{cover} %
\index{cover} (or a minimal right approximation) if in every
factorization \mbox{$\f = \f\psi$}, the map $\mapdef{\psi}{H}{H}$ is
an automorphism. If $\cH$ is closed under extensions, then every
$\cH$-cover is a special precover. This is known as Wakamatsu's lemma;
see Xu~\cite[2.1.1]{xu}.  Dually one defines special
$\cH$-preenvelopes %
\index{preenvelope} and $\cH$-envelopes. %
\index{envelope}

\begin{remark}
  Let $\cx{I}$ be a complex of injective modules as in
  \dfnref{GIm}. Then every differential in $\cx{I}$ is a special
  injective precover of its image; this fact is used in the proof of
  \prpref{count}. Similarly, in a complex $\cx{P}$ of projective
  modules as in \dfnref{GPm}, every differential $\dif[i]{}$ is a
  special projective preenvelope of the cokernel of the previous
  differential $\dif[i+1]{}$.
\end{remark}

\subsubsection*{Totally Reflexive Covers and Envelopes}
\label{ed}

The next result of Avramov and Martsinkovsky~\cite[3.1]{LLAAMr02}
corresponds over a Gorenstein local ring to the existence of maximal
Cohen--Macaulay approximations %
\index{approximation} in the sense of Auslander and
Buchweitz~\cite{MAsROB89}.

\begin{proposition}
  \label{prp:sGprecov}
  Let $R$ be Noetherian. For every finitely generated $R$-module $M$
  of finite G-dimension there is an exact sequence of finitely
  generated modules $0 \to K \to G \to M \to 0,$ where $G$ is totally
  reflexive and one has $\pd{K} = \max\{0,\Gdim{M}-1\}$.  In
  particular, every finitely generated $R$-module of finite
  G-dimension has a special totally reflexive precover.
\end{proposition}

An unpublished result of Auslander states that every finitely
generated module over a Gorenstein local ring has a totally reflexive
cover; see Enochs, Jenda, and Xu~\cite{EJX-99} for a generalization. A
strong converse is contained in the next theorem, which combines
Auslander's result with recent work of several authors; see
\cite{CPST-08} and~\cite{HHlPJr}.

\begin{theorem}
  \label{thm:Gpre}
  For a local ring $\Rmk$ the next conditions are equivalent.
  \begin{description}[$(iii)$]\setlength{\labelsep}{\mylabelsep}
  \item[\hfill $(i)$] Every finitely generated $R$-module has a
    totally reflexive cover.%
    \index{totally reflexive!(pre)cover}
  \item[\hfill $(ii)$] The residue field $k$ has a special totally
    reflexive precover.
  \item[\hfill $(iii)$] Every finitely generated $R$-module has a
    totally reflexive envelope.
  \item[\hfill $(iv)$] Every finitely generated $R$-module has a
    special totally reflexive preenvelope.
  \item[\hfill $(v)$] $R$ is Gorenstein or every totally reflexive
    $R$-module is free.
  \end{description}
\end{theorem}

\noindent
If $R$ is local and Henselian (e.g.\ complete), then existence of a
totally reflexive precover implies existence of a totally reflexive
cover; see \cite[2.5]{RTk05a}. In that case one can drop ``special''
in part $(ii)$ above. In general, though, the next question from
\cite{CPST-08} remains open.

\begin{question}
  Let $\Rmk$ be local. If $k$ has a totally reflexive precover, is
  then $R$ Gorenstein or every totally reflexive $R$-module free?
\end{question}

\subsubsection*{Gorenstein Projective Precovers}

The following result is proved by Holm in~\cite[2.10]{HHl04a}.

\begin{proposition}
  \label{prp:sGPprecov}
  For every $R$-module $M$ of finite Gorenstein projective dimension
  there is an exact sequence $0 \to K \to A \to M \to 0,$ where $A$ is
  Gorenstein projective and $\pd{K} = \max\{0,\Gpd{M}-1\}$.  In
  particular, every $R$-module of finite Gorenstein projective
  dimension has a special Gorenstein projective precover.
\end{proposition}

For an important class of rings, J{\o}rgensen~\cite{PJr07} and Murfet
and Salarian~\cite{DMfSSl} prove existence of Gorenstein projective
precovers for all modules:

\begin{theorem}
  \label{thm:Gprojprecover}
  If $R$ is Noetherian of finite Krull dimension, then every
  $R$-module has a Gorenstein projective precover.%
  \index{Gorenstein projective!precover}
\end{theorem}

\begin{remark}
  Actually, the argument in Krause's proof of \cite[7.12(1)]{HKr05}
  applies to the setup in \cite{PJr07} and yields existence of a
  \emph{special} Gorenstein projective precover for every module over
  a ring as in \thmref{Gprojprecover}.
\end{remark}

Over any ring, every module has a special projective precover; hence:

\begin{prb}
  Describe the rings over which every module has a (special)
  Gorenstein projective precover.
\end{prb}

\subsubsection*{Gorenstein Injective Preenvelopes}

In \cite[2.15]{HHl04a} one finds:

\begin{proposition}
  \label{prp:sGIpreenv}
  For every $R$-module $M$ of finite Gorenstein injective dimension
  there is an exact sequence $0 \to M \to B \to C \to 0,$ where $B$ is
  Gorenstein injective and $\id{C} = \max\{0,\Gid{M}-1\}$.  In
  particular, every $R$-module of finite Gorenstein injective
  dimension has a special Gorenstein injective preenvelope.
\end{proposition}

Over Noetherian rings, existence of Gorenstein injective preenvelopes
for all modules is proved by Enochs and L\'{o}pez-Ramos in
\cite{EEnJLR02}. Krause \cite[7.12]{HKr05} proves a stronger result:

\begin{theorem}
  \label{thm:Gingpreenv}
  \index{Gorenstein injective!preenvelope} If $R$ is Noetherian, then
  every $R$-module has a special Gorenstein injective preenvelope.
\end{theorem}

Over Gorenstein rings, Enochs, Jenda, and Xu \cite[6.1]{EJX-96c} prove
even more:

\begin{proposition}
  If $R$ is Gorenstein of finite Krull dimension, then every
  $R$-module has a Gorenstein injective envelope.
\end{proposition}

Over any ring, every module has an injective envelope; this suggests:

\begin{prb}
  Describe the rings over which every module has a Gorenstein
  injective (pre)envelope.
\end{prb}

Over a Noetherian ring, every module has an injective cover;
see~Enochs \cite[2.1]{EEn81}. A Gorenstein version of this result is
recently established by Holm and J\o rgensen~\cite[3.3(b)]{HHlPJra}:

\begin{proposition}
  If $R$ is a homomorphic image of a Gorenstein ring of finite Krull
  dimension, then every $R$-module has a Gorenstein injective cover.
\end{proposition}

\subsubsection*{Gorenstein Flat Covers}

The following existence result is due to Enochs and
L\'opez-Ramos~\cite[2.11]{EEnJLR02}.

\begin{theorem}
  \label{thm:Gflatcover}
  \index{Gorenstein flat!cover} If $R$ is Noetherian, then every
  $R$-module has a Gorenstein flat cover.
\end{theorem}

\begin{remark}
  \label{rmk:fgprecov}
  Let $R$ be Noetherian and $M$ be a finitely generated $R$-module. If
  $M$ has finite G-dimension, then by \prpref{sGprecov} it has a
  finitely generated Gorenstein projective/flat precover,
  cf.~\prpref{Gfd-gdim}. If $M$ has infinite G-dimension, it still has
  a Gorenstein projective/flat precover by
  \thmref[Theorems~]{Gprojprecover} and \thmref[]{Gflatcover}, but by
  \thmref{Gpre} this need \emph{not} be finitely~generated.
\end{remark}

Over any ring, every module has a flat cover, as proved by Bican, El
Bashir, and Enochs~\cite{BEE-01}. This motivates:

\begin{prb}
  Describe the rings over which every module has a Gorenstein flat
  (pre)cover.
\end{prb}

Over a Noetherian ring, every module has a flat preenvelope;
cf.~Enochs~\cite[5.1]{EEn81}. A Gorenstein version of this result
follows from \thmref[Thm~]{GFprod} and \cite[2.5]{EEnJLR02}:

\begin{proposition}
  If $R$ is a homomorphic image of a Gorenstein ring of finite Krull
  dimension, then every $R$-module has a Gorenstein flat preenvelope.
\end{proposition}

\subsubsection*{Relative Cohomology via Gorenstein Projective Modules}
\label{sec:kd}
\index{relative cohomology}

The notion of a proper resolution is central in relative homological
algebra. Here is a special case:

\begin{dfn}
  An augmented Gorenstein projective resolution, $$\cx{A^+} = \cdots
  \xra{\f_{i+1}} A_i \xra{\f_{i}} A_{i-1} \xra{\f_{i-1}} \cdots
  \xra{\f_{1}} A_0 \xra{\f_{0}} M \lra 0$$ of an $R$-module $M$ is
  said to be \emph{proper} if the complex $\Hom{A'}{\cx{A^+}}$ is
  acyclic for every Gorenstein projective $R$-module $A'$.
\end{dfn}

\begin{remark}
  Assume that every $R$-module has a Gorenstein projective
  precover. Then every $R$-module has a proper Gorenstein projective
  resolution constructed by taking as $\f_0$ a Gorenstein projective
  precover of $M$ and as $\f_i$ a Gorenstein projective precover of
  $\Ker{\f_{i-1}}$ for $i>0$.
\end{remark}

\begin{dfn}
  The \emph{relative Gorenstein projective dimension} %
  \index{Gorenstein projective!relative $\sim$ dimension} of an
  $R$-module \mbox{$M \ne 0$}, denoted by $\rGpd{M}$, is the least
  integer \mbox{$n \geqslant 0$} such that there exists a proper
  Gorenstein projective resolution of $M$ with \mbox{$A_i=0$} for all
  \mbox{$i>n$}.  If no such $n$ (or no such resolution) exists, then
  $\rGpd{M}$ is infinite. By convention, set $\rGpd{0} = -\infty$.
\end{dfn}

The following result is a consequence of \prpref{sGPprecov}.

\begin{proposition}
  \label{prp:rGpd}
  For every $R$-module $M$ one has $\rGpd{M} = \Gpd{M}$.
\end{proposition}

It is shown in \cite[\S 8.2]{rha} that the next definition makes
sense.

\begin{dfn}
  Let $M$ and $N$ be $R$-modules and assume that $M$ has a proper
  Gorenstein projective resolution $\cx{A}$. The $i\,$th relative
  cohomology module $\GPExt{i}{M}{N}$ is $\HH[i]{\Hom{\cx{A}}{N}}$.
\end{dfn}

\begin{remark}
  Let $R$ be Noetherian, and let $M$ and $N$ be finitely generated
  $R$-modules. Unless $M$ has finite G-dimension, it is not clear
  whether the cohomology modules $\GPExt{i}{M}{N}$ are finitely
  generated, cf.~\rmkref{fgprecov}.
\end{remark}

Vanishing of relative cohomology $\operatorname{Ext}^i_{\rm GP}$
characterizes modules of finite Gorenstein projective dimension. The
proof is standard; see \cite[3.9]{HHl08}.

\begin{theorem}
  \label{thm:GPExt}
  Let $M$ be an $R$-module that has a proper Gorenstein projective
  resolution. For every integer $n\ge 0$ the next conditions are
  equivalent.
  \begin{description}[$(iii)$]\setlength{\labelsep}{\mylabelsep}
  \item[\hfill $(i)$] $\Gpd{M} \le n$.
  \item[\hfill $(ii)$] $\GPExt{i}{M}{-}=0$ for all $i>n$.
  \item[\hfill $(iii)$] $\GPExt{n+1}{M}{-}=0$.
  \end{description}
\end{theorem}

\begin{remark}
  Relative cohomology based on totally reflexive modules is studied
  in~\cite{LLAAMr02}. The results that correspond to \prpref{rGpd} and
  \thmref{GPExt} in that setting are contained in \cite[4.8 and
  4.2]{LLAAMr02}.
\end{remark}

\subparagraph{Notes}

\begin{petit}
  Based on a notion of coproper Gorenstein injective resolutions, one
  can define a relative Gorenstein injective dimension and cohomology
  functors $\operatorname{Ext}^i_{\rm GI}$ with properties analogous
  to those of $\operatorname{Ext}^i_{\rm GP}$ described above.  There
  is, similarly, a relative Gorenstein flat dimension and a relative
  homology theory based on proper Gorenstein projective/flat
  resolutions. The relative Gorenstein injective dimension and the
  relative Gorenstein flat dimension were first studied by Enochs and
  Jenda~\cite{EEnOJn94b} for modules over Gorenstein rings. The
  question of balancedness for relative (co)homology is treated by
  Enochs and Jenda~\cite{rha}, Holm~\cite{HHl04b}, and
  Iacob~\cite{AIc07}.
 
  In \cite{LLAAMr02} Avramov and Martsinkovsky also study the
  connection between relative and Tate cohomology for finitely
  generated modules. This study is continued by Veliche \cite{OVl06}
  for arbitrary modules, and a dual theory is developed by Asadollahi
  and Salarian~\cite{JAsSSl06a}. J{\o}rgensen~\cite{PJr07},
  Krause~\cite{HKr05}, and Takahashi~\cite{RTk06b} study connections
  between Gorenstein relative cohomology and generalized notions of
  Tate cohomology.

  Sather-Wagstaff and White~\cite{SSWDWh08} use relative cohomology to
  define an Euler characteristic for modules of finite G-dimension.
  In collaboration with Sharif, they study cohomology theories related
  to generalized Gorenstein dimensions~\cite{SSW-}.
\end{petit}

\section{Modules over Local Homomorphisms}
\label{sec:mlh}

In this section, $\mapdef{\f}{(R,\m)}{(S,\n)}$ is a local
homomorphism, that is, there is a containment $\f(\m) \subseteq
\n$. The topic is Gorenstein dimensions over $R$ of finitely generated
$S$-modules. The utility of this point of view is illustrated by a
generalization, due to Christensen and Iyengar~\cite[4.1]{LWCSIn07},
of the Auslander--Bridger Formula~(\thmref{ABr}):%
\index{Auslander--Bridger Formula}

\begin{theorem}
  Let $N$ be a finitely generated $S$-module. If $N$ has finite
  Gorenstein flat dimension as an $R$-module via $\f$, then there is
  an equality
  $$\Gfd{N} = \dptR - \dpt{N}.$$
\end{theorem}

\noindent
For a finitely generated $S$-module $N$ of finite flat dimension over
$R$, this equality follows from work of Andr\'e \cite[II.57]{hac}.
For a finitely generated $S$-module of finite injective dimension over
$R$, an affirmative answer to the next question is already in
\cite[5.2]{RTkYYs04} by Takahashi and Yoshino.

\begin{question}
  Let $N$ be a non-zero finitely generated $S$-module. If $N$ has
  finite Gorenstein injective dimension as an $R$-module via $\f$,
  does then the equality $\Gid{N} = \dptR$ hold? (For $\f =
  \operatorname{Id}_R$ this is \thmref{Bass_formula}.)
\end{question}

The next result of Christensen and Iyengar \cite[4.8]{LWCSIn07}
should be compared to \thmref{GFcompletion}.

\begin{theorem}
  \label{thm:GFcompletionft}
  Let $N$ be a finitely generated $S$-module. If $N$ has finite
  Gorenstein flat dimension as an $R$-module via $\f$, there is an
  equality $$\Gfd{N} = \Gfd[\Rhat]{\tpp[S]{\Shat}{N}}.$$
\end{theorem}

\subsubsection*{The Frobenius Endomorphism}

If $R$ has positive prime characteristic, we denote by $\phi$ the
Frobenius endomorphism on $R$ and by $\phi^n$ its $n$-fold
composition.  The next two theorems are special cases of \cite[8.14
and 8.15]{SInSSW04} by Iyengar and Sather-Wagstaff and of
\cite[5.5]{HBFAJF07} by Foxby and Frankild.%
\index{Frobenius endomorphism}; together, they constitute the
Gorenstein counterpart of Theorem~F.

\begin{theorem}
  \label{thm:FrobGfd}
  Let $R$ be local of positive prime characteristic. The following
  conditions are equivalent.
  \begin{description}[$(iii)$]\setlength{\labelsep}{\mylabelsep}
  \item[\hfill $(i)$] $R$ is Gorenstein.
  \item[\hfill $(ii)$] $R$ has finite Gorenstein flat dimension as an
    $R$-module via $\phi^n$ for some~\mbox{$n \ge 1$}.
  \item[\hfill $(iii)$] $R$ is Gorenstein flat as an $R$-module via
    $\phi^n$ for every $n \ge 1$.
  \end{description}
\end{theorem}

\begin{theorem}
  \label{thm:FrobGid}
  Let $R$ be local of positive prime characteristic, and assume that
  it is a homomorphic image of a Gorenstein ring.  The following
  conditions are equivalent:
  \begin{description}[$(iii)$]\setlength{\labelsep}{\mylabelsep}
  \item[\hfill $(i)$] $R$ is Gorenstein.%
    \index{Gorenstein ring}
  \item[\hfill $(ii)$] $R$ has finite Gorenstein injective dimension
    as $R$-module via $\phi^n$ for some~\mbox{$n\!\ge\!1$}.
  \item[\hfill $(iii)$] $R$ has Gorenstein injective dimension equal
    to $\dimR$ as an $R$-module via $\phi^n$ for every $n \ge 1$.
  \end{description}
\end{theorem}

\begin{petit}
  \noindent
  Part $(iii)$ in \thmref{FrobGid} is actually not included in
  \cite[5.5]{HBFAJF07}. Part $(iii)$ follows, though, from part $(i)$
  by \corref{Gor-Gid} and \thmref{Gid-Chouinard}.
\end{petit}

\subsubsection*{G-dimension over a Local Homomorphism}

The homomorphism $\mapdef{\f}{(R,\m)}{(S,\n)}$ fits in a commutative
diagram of local homomorphisms:
$$\xymatrix@=1.5em{{} & R' \ar@{->>}[dr]^-{\varphi'} & {} \\
  R \ar[ur]^-{\dot{\varphi}} \ar[r]_-{\varphi} & S\;
  \ar@{^(->}[r]_-{\iota} & \Shat}$$ where $\dot{\f}$ is flat with
regular closed fiber $R'/\m R'$, the ring $R'$ is complete, and $\f'$
is surjective. Set $\grave{\f}=\iota\f$; a diagram as above is called
a \emph{Cohen factorization of} $\grave{\f}$. This is a construction
due to Avramov, Foxby, and Herzog~\cite[1.1]{AFH-94}.%
\index{Cohen factorization}

The next definition is due to Iyengar and Sather-Wagstaff~\cite[\S
3]{SInSSW04}; it is proved \emph{ibid.} that it is independent of the
choice of Cohen factorization.
\begin{dfn}
  \label{dfn:gdimf}
  Choose a Cohen factorization of $\grave{\f}$ as above. For a
  finitely generated $S$-module $N$, the \emph{G-dimension of $N$ over
    $\f$} %
  \index{G-dimension!over a local homomorphism} is given as
$$\Gdim[\f]{N} = \Gdim[R']{\tpp[S]{\Shat}{N}} - \edim{(R'/\m R')}.$$
\end{dfn}

\begin{example}
  \label{exa:k}
  Let $k$ be a field and let $\f$ be the extension from $k$ to the
  power series ring $\pows{x}$. Then one has $\Gdim[\f]{\pows{x}} =
  -1$.
\end{example}

Iyengar and Sather-Wagstaff~\cite[8.2]{SInSSW04} prove:

\begin{theorem}
  \label{thm:gdim-Gfd}
  Assume that $R$ is a homomorphic image of a Gorenstein ring. A
  finitely generated $S$-module $N$ has finite G-dimension over $\f$
  if and only if it has finite Gorenstein flat dimension as an
  $R$-module via~$\f$.
\end{theorem}

\noindent
It is clear from \exaref{k} that $\Gfd{N}$ and $\Gdim[\f]{N}$ need not
be equal.

\section{Local Homomorphisms of Finite G-dimension}
\label{sec:lrhgd}

This section treats transfer of ring theoretic properties along a
local homomorphism of finite G-dimension. Our focus is on the
Gorenstein property, which was studied by Avramov and Foxby in
\cite{LLAHBF97}, and the Cohen--Macaulay property, studied by Frankild
in \cite{AFr01}.

As in \secref{mlh}, $\mapdef{\f}{(R,\m)}{(S,\n)}$ is a local
homomorphism. In view of \dfnref{gdimf}, a notion from
\cite[4.3]{LLAHBF97} can be defined as follows:

\begin{dfn}
  Set $\Gdim[]{\f}=\Gdim[\f]{S}$; the homomorphism $\f$ is said to be
  of \emph{finite G-dimension} %
  \index{local homomorphism!of finite G-dimension} if this number is
  finite.
\end{dfn}

\begin{remark}
  The homomorphism $\f$ has finite G-dimension if $S$ has finite
  Gorenstein flat dimension as an $R$-module via $\f$, and the
  converse holds if $R$ is a homomorphic image of a Gorenstein ring.
  This follows from \thmref[Theorems~]{GFcompletionft} and
  \thmref[]{gdim-Gfd}, in view of \cite[3.4.1]{SInSSW04},
\end{remark}

The next descent result is \cite[4.6]{LLAHBF97}.

\begin{theorem}
  \label{thm:GorAD}
  Let $\f$ be of finite G-dimension, and assume that $R$ is a
  homomorphic image of a Gorenstein ring. For every $S$-module $N$ one~has:
  \begin{description}[$(a)$]\setlength{\labelsep}{\mylabelsep}
  \item[\hfill \rm (a)] If\, $\fd[S]{N}$ is finite then $\Gfd{N}$ is
    finite.
  \item[\hfill \rm (b)] If\, $\id[S]{N}$ is finite then $\Gid{N}$ is
    finite.
  \end{description}
\end{theorem}

It is not known if the composition of two local homomorphisms of
finite G-dimension has finite G-dimension, but it would follow from an
affirmative answer to \qstref{trans}, cf.~\cite[4.8]{LLAHBF97}.%
\index{G-dimension!Transitivity Question} Some
insight is provided by \thmref{qGorcompos} and the next result, which
is due to Iyengar and Sather-Wagstaff \cite[5.2]{SInSSW04}.

\begin{theorem}
  \label{thm:compos}
  Let $\mapdef{\psi}{S}{T}$ be a local homomorphism such that
  $\fd[S]{T}$ is finite. Then $\Gdim[]{\psi\f}$ is finite if and only
  if $\Gdim[]{\f}$ is finite.
\end{theorem}

\subsubsection*{Quasi-Gorenstein Homomorphisms}
\label{sec:qgh}

Let $M$ be a finitely generated module over a local ring $\Rmk$. For
every integer $i \ge 0$ the $i\,$th \emph{Bass number} $\bas{i}{M}$ is
the dimension of the $k$-vector space $\Ext{i}{k}{M}$.

\begin{dfn}
  \label{def:QGorHomo}
  The homomorphism $\f$ is called \emph{quasi-Gorenstein} %
  \index{local homomorphism!quasi-Gorenstein} if it has finite
  G-dimension and for every $i \ge 0$ there is an equality of Bass
  numbers $$\bas{i+\dptR}{R} = \bas[S]{i+\dpt[]{S}}{S}.$$
\end{dfn}

\begin{example}
  \label{exa:qgor}
  If $R$ is Gorenstein, then the natural surjection $R
  \twoheadrightarrow k$ is quasi-Gorenstein.  More generally, if $\f$
  is surjective and $S$ is quasi-perfect%
  \index{quasi-perfect module} as an $R$-module via $\f$, then $\f$ is
  quasi-Gorenstein if and only if there is an isomorphism
  $\Ext{g}{S}{R} \is S$ where $g = \Gdim{S}$ with; see \cite[6.5, 7.1,
  7.4]{LLAHBF97}.
\end{example}

Several characterizations of the quasi-Gorenstein property are given
in \cite[7.4 and 7.5]{LLAHBF97}. For example, it is sufficient that
$\Gdim[]{\f}$ is finite and the equality of Bass numbers holds for
some $i >0$.

The next ascent-descent result is~\mbox{\cite[7.9]{LLAHBF97}.}

\begin{theorem}
  \label{thm:qGorAD}
  Let $\f$ be quasi-Gorenstein and assume that $R$ is a homomorphic
  image of a Gorenstein ring. For every $S$-module $N$ one has:
  \begin{description}[$(a)$]\setlength{\labelsep}{\mylabelsep}
  \item[\hfill\rm (a)] $\Gfd[S]{N}$ is finite if and only if $\Gfd{N}$
    is finite.
  \item[\hfill\rm (b)] $\Gid[S]{N}$ is finite if and only if $\Gid{N}$
    is finite.
  \end{description}
\end{theorem}

Ascent and descent of the Gorenstein property is described by
\cite[7.7.2]{LLAHBF97}. It should be compared to part (b) in
Theorem~G.

\begin{theorem}
  \label{thm:GorAscDescFtGFD}
  The following conditions on $\f$ are equivalent.
  \begin{description}[$(iii)$]\setlength{\labelsep}{\mylabelsep}
  \item[\hfill $(i)$] $R$ and $S$ are Gorenstein.%
    \index{Gorenstein ring!local}%
  \item[\hfill $(ii)$] $R$ is Gorenstein and $\f$ is quasi-Gorenstein.
  \item[\hfill $(iii)$] $S$ is Gorenstein and $\f$ is of finite
    G-dimension.
  \end{description}
\end{theorem}

The following (de)composition result is \cite[7.10, 8.9, and
8.10]{LLAHBF97}. It should be compared to \thmref{trans}.

\begin{theorem}
  \label{thm:qGorcompos}
  Assume that $\f$ is quasi-Gorenstein, and let $\mapdef{\psi}{S}{T}$
  be a local homomorphism. The following assertions hold.
  \begin{description}[$(a)$]\setlength{\labelsep}{\mylabelsep}
  \item[\hfill \rm (a)] $\Gdim[]{\psi\f}$ is finite if and only if
    $\Gdim[]{\psi}$ is finite.
  \item[\hfill \rm (b)] $\psi\f$ is quasi-Gorenstein if and only if
    $\psi$ is quasi-Gorenstein.
  \end{description}
\end{theorem}

\subsubsection*{Quasi-Cohen--Macaulay Homomorphisms}
\label{subsec:CMAscDescThms}

The next definition from \cite[5.8 and 6.2]{AFr01} uses terminology
from \dfnref{qperf} and the remarks before \dfnref{gdimf}.

\begin{dfn}
  The homomorphism $\f$ is \emph{quasi-Cohen--Macaulay}, %
  \index{local homomorphism!quasi-Cohen--Macaulay} for short
  \emph{quasi-CM}, if $\grave{\f}$ has a Cohen factorization where
  $\Shat$ is quasi-perfect over $R'$.
\end{dfn}

\noindent
If $\f$ is quasi-CM, then $\Shat$ is a quasi-perfect $R'$-module in
every Cohen factorization of $\grave{\f}$; see \cite[5.8]{AFr01}. The
following theorem is part of \cite[6.7]{AFr01}; it should be compared
to part (a) in Theorem~G.

\begin{theorem}
  \label{thm:CMAscDsc}
  The following assertions hold.
  \begin{description}[$(a)$]\setlength{\labelsep}{\mylabelsep}
  \item[\hfill\rm (a)] If $R$ is Cohen--Macaulay and $\f$ is quasi-CM,
    then $S$ is Cohen--Macaulay.
  \item[\hfill\rm (b)] If $S$ is Cohen--Macaulay and $\Gdim[]{\f}$ is
    finite, then $\f$ is quasi-CM.
  \end{description}
\end{theorem}

In view of \thmref{qGorcompos}, Frankild's work \cite[6.4 and
6.5]{AFr01} yields:

\begin{theorem}
  Assume that $\f$ is quasi-Gorenstein, and let $\mapdef{\psi}{S}{T}$
  be a local homomorphism. Then $\psi\f$ is quasi-CM if and only if
  $\psi$ is quasi-CM.
\end{theorem}

\subparagraph{Notes}

The composition question addressed in the remarks before
\thmref{compos} is investigated further by Sather-Wagstaff
\cite{SSW08}.


\section{Reflexivity and Finite G-dimension}
\label{sec:fingdim}

In this section $R$ is Noetherian. Let $M$ be a finitely generated
$R$-module. If $M$ is totally reflexive, then the cohomology modules
$\Ext{i}{M}{R}$ vanish for all $i >0$. The converse is true if $M$ is
known \emph{a priori} to have finite G-dimension,
cf.~\corref{gdim}. In general, though, one can not infer from such
vanishing that $M$ is totally reflexive---explicit examples to this
effect are constructed by Jorgensen and \c{S}ega in
\cite{DAJLMS06}---and this has motivated a search for alternative
criteria for finiteness of G-dimension.

\subsubsection*{Reflexive Complexes}
\label{refl}

One such criterion was given by Foxby and published in
\cite{SYs95b}. Its habitat is the derived category $\Du$ of the
category of $R$-modules. The objects in $\Du$ are $R$-complexes, and
there is a canonical functor $\mathrm{F}$ from the category of
$R$-complexes to $\Du$. This functor is the identity on objects and it
maps homology isomorphisms to isomorphisms in $\Du$. The restriction
of $\mathrm{F}$ to modules is a full embedding of the module category
into $\Du$.%
\index{complex}
\index{derived category}

The homology $\H{\cx{M}}$ of an $R$-complex $\cx{M}$ is a (graded)
$R$-module, and $\cx{M}$ is said to have \emph{finitely generated
  homology} if this module is finitely generated. That is, if every
homology module $\H[i]{\cx{M}}$ is finitely generated and only
finitely many of them are non-zero.

For $R$-modules $M$ and $N$, the (co)homology of the derived Hom and
tensor product complexes gives the classical Ext and Tor
modules: $$\Ext{i}{M}{N} = \HH[i]{\RHom{M}{N}} \quad\text{and}\quad
\Tor{i}{M}{N} = \H[i]{\Ltp{M}{N}}.$$

\begin{dfn}
  An $R$-complex $\cx{M}$ is \emph{reflexive} %
  \index{reflexive complex} if $\cx{M}$ and $\RHom{\cx{M}}{R}$ have
  finitely generated homology and the canonical morphism $$\cx{M} \lra
  \RHom{\RHom{\cx{M}}{R}}{R}$$ is an isomorphism in the derived
  category $\Du$.  The full subcategory of $\Du$ whose objects are the
  reflexive $R$-complexes is denoted by $\R$.
\end{dfn}

\begin{petit}
  \noindent
  The requirement in the definition that the complex
  $\RHom{\cx{M}}{R}$ has finitely generated homology is redundant but
  retained for historical reasons; see \cite[3.3]{AIL-}.
\end{petit}

\enlargethispage*{\baselineskip}

\thmref{Rmod} below is Foxby's criterion for finiteness of G-dimension
of a finitely generated module~\cite[2.7]{SYs95b}. It differs
significantly from \dfnref{G-dim} as it does not involve construction
of a G-resolution of the module.

\begin{theorem}
  \label{thm:Rmod}
  Let $R$ be Noetherian. A finitely generated $R$-module has finite
  G-dimension if and only if it belongs to $\R$.
\end{theorem}

If $R$ is local, then the next result is \cite[2.3.14]{lnm}. In the
generality stated below it follows from \thmref[Theorems~]{Goto} and
\thmref[]{Rmod}: the implication $(ii)\Rightarrow (iii)$ is the least
obvious, it uses \cite[2.1.12]{lnm}.

\begin{corollary}
  Let $R$ be Noetherian. The following conditions are equivalent.
  \begin{description}[$(iii)$]\setlength{\labelsep}{\mylabelsep}
  \item[\hfill $(i)$] $R$ is Gorenstein.
  \item[\hfill $(ii)$] Every $R$-module is in $\R$.
  \item[\hfill $(iii)$] Every $R$-complex with finitely generated
    homology is in $\R$.
  \end{description}
\end{corollary}

\subsubsection*{G-dimension of Complexes}
\label{gdimc}

Having made the passage to the derived category, it is natural to
consider G-dimension for complexes. For every $R$-complex $\cx{M}$
with finitely generated homology there exists a complex $\cx{G}$ of
finitely generated free $R$-modules, which is isomorphic to $\cx{M}$
in $\Du$; see \cite[1.7(1)]{LLAHBF91}. In Christensen's
\cite[ch.~2]{lnm} one finds the next definition and the two theorems
that follow.

\begin{dfn}
  Let $\cx{M}$ be an $R$-complex with finitely generated homology. If
  $\H{\cx{M}}\ne 0$, then the \emph{G-dimension} %
  \index{G-dimension} of $\cx{M}$ is the least integer $n$ such that
  there exists a complex $\cx{G}$ of totally reflexive $R$-modules
  which is isomorphic to $\cx{M}$ in $\Du$ and has $G_i=0$ for all $i
  > n$. If no such integer $n$ exists, then $\Gdim{\cx{M}}$ is
  infinite. If $\H{\cx{M}}= 0$, then $\Gdim{\cx{M}} = -\infty$ by
  convention.
\end{dfn}

\noindent Note that this extends \dfnref{G-dim}. As a common
generalization of \thmref{Rmod} and \corref{gdim} one has
\cite[2.3.8]{lnm}:

\begin{theorem}
  \label{thm:R}
  Let $R$ be Noetherian. An $R$-complex $\cx{M}$ with finitely
  generated homology has finite G-dimension if an only if it is
  reflexive.  Furthermore, for every reflexive $R$-complex $\cx{M}$
  there is an equality $$\Gdim{\cx{M}} =
  \sup\setof{i\in\ZZ}{\HH[i]{\RHom{\cx{M}}{R}}\ne 0}.$$
\end{theorem}

In broad terms, the theory of G-dimension for finitely generated
modules extends to complexes with finitely generated homology.  One
example is the next extension \cite[2.3.13]{lnm} of the
Auslander--Bridger Formula~(\thmref{ABr}). %
\index{Auslander--Bridger Formula}

\begin{theorem}
  Let $R$ be local. For every complex $\cx{M}$ in $\R$ one
  has $$\Gdim{\cx{M}} = \dptR - \dpt{\cx{M}}.$$
\end{theorem}

\noindent
Here the \emph{depth} %
\index{depth} of a complex $\cx{M}$ over a local ring $\Rmk$ is
defined by extension of the definition for modules, that
is, $$\dpt{\cx{M}} = \inf\setof{i\in\ZZ}{\HH[i]{\RHom{k}{\cx{M}}}\ne
  0}.$$

\subparagraph{Notes}

\begin{petit}
  In \cite[ch.~2]{lnm} the theory of G-dimension for complexes is
  developed in detail.

  Generalized notions of G-dimension---based on reflexivity with
  respect to semidualizing modules and complexes---are studied by
  Avramov, Iyengar, and Lipman~\cite{AIL-}, Christensen~\cite{LWC01a},
  Frankild and Sather-Wagstaff~\cite{AFrSSW07a}, Gerko~\cite{AAG01a},
  Golod~\cite{ESG84}, Holm and J{\o}rgensen~\cite{HHlPJr07a}, by
  Salarian, Sather-Wagstaff, and Yassemi~\cite{SSY-06}, and
  White~\cite{DWh}. See also the notes in \secref{gdim}.
\end{petit}

\section{Detecting Finiteness of Gorenstein Dimensions}
\label{sec:fingd}

In the previous section, we discussed a resolution-free
characterization of modules of finite G-dimension (\thmref{Rmod}). The
topic of this section is similar characterizations of modules of
finite Gorenstein projective/injective/flat dimension. By work of
Christensen, Frankild, and Holm~\cite{CFH-06}, appropriate criteria
are available for modules over a Noetherian ring that has a dualizing
complex (\thmref[Theorems~]{Amod} and \thmref[]{Bmod}). As mentioned
in the Introduction, a Noetherian ring has a dualizing complex if and
only if it is a homomorphic image of a Gorenstein ring of finite Krull
dimension. For example, every complete local ring has a dualizing
complex by Cohen's structure theorem.

\subsubsection*{Auslander Categories}
\label{ab}

The next definition is due to Foxby; see~\cite[3.1]{LLAHBF97} and
\cite[\S 2]{HBF72}.

\begin{dfn}
  \label{dfn:A}
  Let $R$ be Noetherian and assume that it has a dualizing complex %
  \index{dualizing complex} $\cx{D}$. The \emph{Auslander class} %
  \index{Auslander class} $\A$ is the full subcategory of the derived
  category $\Du$ whose objects $\cx{M}$ satisfy the following
  conditions.
  \begin{description}[$(3)$]\setlength{\labelsep}{\mylabelsep}
  \item[\hfill $(1)$] $\H[i]{\cx{M}} =0$ for $|i| \gg 0$.
  \item[\hfill $(2)$] $\H[i]{\Ltp{\cx{D}}{\cx{M}}}=0$ for $i \gg 0$.
  \item[\hfill $(3)$] The natural map $\cx{M}\to
    \RHom{\cx{D}}{\Ltp{\cx{D}}{\cx{M}}}$ is invertible in~$\Du$.
  \end{description}
\end{dfn}

The relation to Gorenstein dimensions is given by \cite[4.1]{CFH-06}:

\begin{theorem}
  \label{thm:Amod}
  Let $R$ be Noetherian and assume that it has a dualizing
  complex. For every $R$-module $M$, the following conditions are
  equivalent.
  \begin{description}[$(iii)$]\setlength{\labelsep}{\mylabelsep}
  \item[\hfill $(i)$] $M$ has finite Gorenstein projective dimension.%
    \index{Gorenstein projective!dimension}
  \item[\hfill $(ii)$] $M$ has finite Gorenstein flat dimension.
    \index{Gorenstein flat!dimension}
  \item[\hfill $(iii)$] $M$ belongs to $\A$.
  \end{description}
\end{theorem}

\begin{remark}
  \label{rmk:resfree1}
  The equivalence of $(i)/(ii)$ and $(iii)$ in \thmref{Amod} provides
  a resolution-free characterization of modules of finite Gorenstein
  projective/flat dimension over a ring that has a dualizing
  complex. Every complete local ring has a dualizing complex, so in
  view of \thmref{GPcompletion}/\thmref[]{GFcompletion} there is
   a resolution-free characterization of modules of finite
  Gorenstein projective/flat dimension over any local ring.
\end{remark}

The next definition is in \cite[3.1]{LLAHBF97}; the theorem that
follows is \cite[4.4]{CFH-06}.

\begin{dfn}
  \label{dfn:B}
  Let $R$ be Noetherian and assume that it has a dualizing complex
  $\cx{D}$. The \emph{Bass class} %
  \index{Bass class} $\B$ is the full subcategory of the derived
  category $\Du$ whose objects $\cx{M}$ satisfy the following
  conditions.
  \begin{description}[$(3)$]\setlength{\labelsep}{\mylabelsep}
  \item[\hfill $(1)$] $\HH[i]{\cx{M}} =0$ for $|i| \gg 0$.
  \item[\hfill $(2)$] $\HH[i]{\RHom{\cx{D}}{\cx{M}}}=0$ for $i \gg 0$.
  \item[\hfill $(3)$] The natural map
    $\Ltp{\cx{D}}{\RHom{\cx{D}}{\cx{M}}}\to\cx{M}$ is invertible
    in~$\Du$.
  \end{description}
\end{dfn}

\begin{theorem}
  \label{thm:Bmod}
  Let $R$ be Noetherian and assume that it has a dualizing
  complex. For every $R$-module $M$, the following conditions are
  equivalent.
  \begin{description}[$(ii)$]\setlength{\labelsep}{\mylabelsep}
  \item[\hfill $(i)$] $M$ has finite Gorenstein injective dimension.%
    \index{Gorenstein injective!dimension}
  \item[\hfill $(ii)$] $M$ belongs to $\B$.
  \end{description}
\end{theorem}

From the two theorems above and from
\thmref[Theorems~]{Iwanaga-Gorenstein} and
\thmref[]{GF-Iwanaga-Gorenstein} one gets:

\begin{corollary}
  Let $R$ be Noetherian and assume that it has a dualizing
  complex. The following conditions are equivalent.
  \begin{description}[$(iii)$]\setlength{\labelsep}{\mylabelsep}
  \item[\hfill $(i)$] $R$ is Gorenstein.%
    \index{Gorenstein ring}
  \item[\hfill $(ii)$] Every $R$-complex $\cx{M}$ with $\H[i]{\cx{M}}
    =0$ for $|i| \gg 0$ belongs to $\A$.
  \item[\hfill $(iii)$] Every $R$-complex $\cx{M}$ with $\H[i]{\cx{M}}
    =0$ for $|i| \gg 0$ belongs to $\B$.
  \end{description}
\end{corollary}

\subsubsection*{Gorenstein Dimensions of Complexes}
\label{gdc}

It turns out to be convenient to extend the Gorenstein dimensions to
complexes; this is illustrated by \thmref{GidCoBase} below.

In the following we use the notion of a semi-projective
resolution. Every complex has such a resolution, by
\cite[1.6]{LLAHBF91}, and a projective resolution of a module is
semi-projective. In view of this and \thmref{Gpd}, the next
definition, which is due to Veliche~\cite[3.1 and 3.4]{OVl06}, extends
\dfnref{Gpd}.

\begin{dfn}
  Let $\cx{M}$ be an $R$-complex. If $\H{\cx{M}}\ne 0$, then the
  \emph{Gorenstein projective dimension} %
  \index{Gorenstein projective!dimension} of $\cx{M}$ is the least
  integer $n$ such that $\H[i]{\cx{M}}=0$ for all $i > n$ and there
  exists a semi-projective resolution $\cx{P}$ of $\cx{M}$ for which
  the module $\Coker{(P_{n+1} \to P_n)}$ is Gorenstein projective. If
  no such $n$ exists, then $\Gpd{\cx{M}}$ is infinite. If $\H{\cx{M}}=
  0$, then $\Gpd{\cx{M}} = -\infty$ by convention.
\end{dfn}

In the next theorem, which is due to Iyengar and
Krause~\cite[8.1]{SInHKr06}, the \emph{unbounded Auslander class} %
\index{Auslander class} $\smash{\Au}$ is the full subcategory of $\Du$
whose objects satisfy conditions $(2)$ and $(3)$ in \dfnref{A}.

\begin{theorem}
  \label{thm:Acx}
  Let $R$ be Noetherian and assume that it has a dualizing
  complex. For every $R$-complex $\cx{M}$, the following conditions
  are equivalent.
  \begin{description}[$(ii)$]\setlength{\labelsep}{\mylabelsep}
  \item[\hfill $(i)$] $\cx{M}$ has finite Gorenstein projective
    dimension.
  \item[\hfill $(ii)$] $\cx{M}$ belongs to $\Au$.
  \end{description}
\end{theorem}

One finds the next definition in \cite[2.2 and 2.3]{JAsSSl06a} by
Asadollahi and Salarian. It uses the notion of a semi-injective
resolution. Every complex has such a resolution, by
\cite[1.6]{LLAHBF91}, and an injective resolution of a module is
semi-injective. In view of \thmref{Gid}, the following extends
\dfnref{Gid}.

\begin{dfn}
  Let $\cx{M}$ be an $R$-complex. If $\H{\cx{M}}\ne 0$, then the
  \emph{Gorenstein injective dimension} %
  \index{Gorenstein injective!dimension} of $\cx{M}$ is the least
  integer $n$ such that $\HH[i]{\cx{M}}=0$ for all $i > n$ and there
  exists a semi-injective resolution $\cx{I}$ of $\cx{M}$ for which
  the module $\Ker{(I^{n} \to I^{n+1})}$ is Gorenstein injective. If
  no such integer $n$ exists, then $\Gid{\cx{M}}$ is infinite. If
  $\H{\cx{M}}= 0$, then $\Gid{\cx{M}} = -\infty$ by convention.
\end{dfn}

In the next result, which is \cite[8.2]{SInHKr06}, the \emph{unbounded
  Bass class} %
\index{Bass class} $\smash{\Bu}$ is the full subcategory of $\Du$
whose objects satisfy  $(2)$ and $(3)$ in \dfnref{B}.

\begin{theorem}
  \label{thm:Bcx}
  Let $R$ be Noetherian and assume that it has a dualizing
  complex. For every $R$-complex $\cx{M}$, the following conditions
  are equivalent.
  \begin{description}[$(ii)$]\setlength{\labelsep}{\mylabelsep}
  \item[\hfill $(i)$] $\cx{M}$ has finite Gorenstein injective
    dimension.
  \item[\hfill $(ii)$] $\cx{M}$ belongs to $\Bu$.
  \end{description}
\end{theorem}

The next result is~\cite[1.7]{LWCSSWc}; it should be compared to
\thmref{GI_cotorsion}.

\begin{theorem}
  \label{thm:GidCoBase}
  Let $R$ be local. For every $R$-module $M$ there is an
  equality $$\Gid{M} = \Gid[\Rhat]{\RHom{\Rhat}{M}}.$$
\end{theorem}

\begin{remark}
  \label{rmk:Gid}
  Via this result, \thmref{Bcx} gives a resolution-free
  characterization of modules of finite Gorenstein injective dimension
  over any local~ring. 
\end{remark}

\subsubsection*{Acyclicity Versus Total Acyclicity}
\label{actac}

The next results characterize Gorenstein rings in terms of the
complexes that define Gorenstein projective/injective/flat
modules. The first one is~\cite[5.5]{SInHKr06}.

\begin{theorem}
  \index{acyclic complex} %
  Let $R$ be Noetherian and assume that it has a dualizing
  complex. Then the following conditions are equivalent.
  \begin{description}[$(iii)$]\setlength{\labelsep}{\mylabelsep}
  \item[\hfill $(i)$] $R$ is Gorenstein.%
    \index{Gorenstein ring}
  \item[\hfill $(ii)$] For every acyclic complex $\cx{P}$ of
    projective $R$-modules and every projective $R$-module $Q$, the
    complex $\Hom{\cx{P}}{Q}$ is acyclic.
  \item[\hfill $(iii)$] For every acyclic complex $\cx{I}$ of
    injective $R$-modules and every injective $R$-module $E$, the
    complex $\Hom{E}{\cx{I}}$ is acyclic.
  \end{description}
\end{theorem}

\noindent
In the terminology of \cite{SInHKr06}, part $(ii)/(iii)$ above says
that every acyclic complex of projective/injective modules is totally
acyclic.%
\index{totally acyclic complex}

The final result is due to Christensen and Veliche~\cite{LWCOVl08}:

\begin{theorem}
  Let $R$ be Noetherian and assume that it has a dualizing
  complex. Then there exist acyclic complexes $\cx{F}$ and $\cx{I}$ of
  flat $R$-modules and injective $R$-modules, respectively, such that
  the following conditions are equivalent.
  \begin{description}[$(iii)$]\setlength{\labelsep}{\mylabelsep}
  \item[\hfill $(i)$] $R$ is Gorenstein.
  \item[\hfill $(ii)$] For every injective $R$-module $E$, the complex
    $\tp{E}{\cx{F}}$ is acyclic.
  \item[\hfill $(iii)$] For every injective $R$-module $E$, the
    complex $\Hom{E}{\cx{I}}$ is acyclic.
  \end{description}
\end{theorem}

\noindent
The complexes $\cx{F}$ and $\cx{I}$ in the theorem have explicit
constructions. It is not known, in general, if there is an explicit
construction of an acyclic complex $\cx{P}$ of projective $R$-modules
such that $R$ is Gorenstein if $\Hom{\cx{P}}{Q}$ is acyclic for every
projective $R$-module $Q$.

\subparagraph{Notes}

\begin{petit}
  In broad terms, the theory of Gorenstein dimensions for modules
  extends to complexes. It is developed in detail by Asadollahi and
  Salarian~\cite{JAsSSl06a}, Christensen, Frankild, and
  Holm~\cite{CFH-06} and \cite{LWCHHl09a}, Christensen and
  Sather-Wagstaff~\cite{LWCSSWc}, and by Veliche~\cite{OVl06}.

  Objects in Auslander and Bass classes with respect to semi-dualizing
  complexes have interpretations in terms of generalized Gorenstein
  dimensions; see \cite{EEnOJn04} and \cite{EEnOJn05} by Enochs and
  Jenda, \cite{HHlPJr06} by Holm and J{\o}rgensen, and \cite{SSW08} by
  Sather-Wagstaff.

  Sharif, Sather-Wagstaff, and White~\cite{SSW-08} study totally
  acyclic complexes of Gorenstein projective modules. They show that
  the cokernels of the differentials in such complexes are Gorenstein
  projective. That is, a ``Gorenstein Gorenstein projective'' module
  is Gorenstein projective.
\end{petit}

\bibliographystyle{spmpsci}

\def\cprime{$'$}
\newcommand{\arxiv}[2][AC]{\href{http://arxiv.org/abs/#2}{\sf
      arXiv:#2 [math.#1]}}
\newcommand{\oldarxiv}[2][AC]{\href{http://arxiv.org/abs/math/#2}{\sf
      arXiv:math/#2
      [math.#1]}}\providecommand{\MR}[1]{\mbox{\href{http://www.ams.org/mathscine%
      t-getitem?mr=#1}{#1}}}
\renewcommand{\MR}[1]{\mbox{\href{http://www.ams.org/mathscinet-getitem?mr=#%
      1}{#1}}}

\printindex
\end{document}